\documentclass[fontsize=12pt,a4paper,headings=normal,
twoside=false,leqno,parskip=half-,abstract=true]{scrartcl}
\usepackage[english]{babel}
\usepackage[utf8]{inputenc}
\usepackage{hyperref}
\hypersetup{
 pdftitle={meanders with three noses},
 pdfauthor={Bernold Fiedler},
 colorlinks=true,
 linkcolor=blue,
 citecolor=blue,
 filecolor=blue,
 urlcolor=blue}
  
 \usepackage{afterpage}

\usepackage{graphicx}
\usepackage{wrapfig} 
\usepackage[format=plain,labelfont=bf,font=small]{caption}
\usepackage{xcolor}
\usepackage[arrow, matrix, curve]{xy}
\usepackage{float}

\usepackage{caption}
\usepackage{tabulary}
\usepackage{array}
\newcolumntype{N}[1]{>{\centering\arraybackslash}m{#1}}

\usepackage{amsmath,amsthm}
\swapnumbers 
\usepackage{amssymb} 

\makeatletter
\newcommand{\tpitchfork}{%
  \vbox{
    \baselineskip\z@skip
    \lineskip-.52ex
    \lineskiplimit\maxdimen
    \m@th
    \ialign{##\crcr\hidewidth\smash{$-$}\hidewidth\crcr$\pitchfork$\crcr}
  }%
}
\makeatother
\usepackage{latexsym}
\usepackage{enumerate}

\usepackage[notref,notcite,color, final 
]{showkeys}

\definecolor{refkey}{rgb}{1,0,0}
\definecolor{labelkey}{rgb}{1,0,0}

\usepackage{cancel}

\usepackage{tikz}

  \mathchardef\ordinarycolon\mathcode`\:
  \mathcode`\:=\string"8000
  \begingroup \catcode`\:=\active
    \gdef:{\mathrel{\mathop\ordinarycolon}}
  \endgroup

\theoremstyle{plain}
\newtheorem{thm}{Theorem}[section]
\newtheorem{lem}[thm]{Lemma}
\newtheorem{prop}[thm]{Proposition}
\newtheorem{cor}[thm]{Corollary}
\newtheorem{defi}[thm]{Definition}

\begin{document}

\title{\LARGE{Design of Sturm global attractors 2:\\
Time-reversible Chafee-Infante lattices\\ 
of 3-nose meanders}
\vspace{1cm}}
{\subtitle{	
	\vspace{1ex}
	{}}\vspace{1ex}
	}

\author{
 \\
\emph{-- Dedicated to Giorgio Fusco, mentor and friend  --}\\
{~}\\
Bernold Fiedler* and Carlos Rocha{**}\\
\vspace{2cm}}

\date{\small{version of \today}}
\maketitle
\thispagestyle{empty}

\vfill

*\\
Institut für Mathematik\\
Freie Universität Berlin\\
Arnimallee 3\\ 
14195 Berlin, Germany\\
 \\
{**}\\
Instituto Superior T\'ecnico\\
Avenida Rovisco Pais\\ 
1049-001 Lisboa, Portugal


\newpage
\pagestyle{plain}
\pagenumbering{roman}
\setcounter{page}{1}

\begin{abstract}
\noindent
This sequel continues our exploration \cite{firo23} of a deceptively ``simple'' class of global attractors, called \emph{Sturm} due to nodal properties.
They arise for the semilinear scalar parabolic PDE
	\begin{equation}\label{eq:*}
	u_t = u_{xx} + f(x,u,u_x) \tag{$*$}
	\end{equation}
on the unit interval $0 < x<1$, under Neumann boundary conditions.
This models the interplay of reaction, advection, and diffusion.

\smallskip\noindent
Our classification is based on the Sturm meanders, which arise from a shooting approach to the ODE boundary value problem of equilibrium solutions $u=v(x)$.
Specifically, we address meanders with only three ``noses'', each of which is innermost to a nested family of upper or lower meander arcs.
The Chafee-Infante paradigm of 1974, with cubic nonlinearity $f=f(u)$, features just two noses.

\smallskip\noindent
We present, and fully prove, a precise description of global PDE connection graphs, graded by Morse index, for such gradient-like Morse-Smale systems \eqref{eq:*}.
The directed edges denote PDE heteroclinic orbits $v_1 \leadsto v_2$ between equilibrium vertices $v_1, v_2$ of adjacent Morse index.
The connection graphs can be described as a lattice-like structure of Chafee-Infante subgraphs.
However, this simple description requires us to adjoin a single ``equilibrium'' vertex, formally, at Morse level -1.
Surprisingly, for parabolic PDEs based on irreversible diffusion, the connection graphs then also exhibit global time reversibility.

\end{abstract}

\tableofcontents


\newpage
\pagenumbering{arabic}
\setcounter{page}{1}

\section{Introduction and main results} \label{Intro}

\numberwithin{equation}{section}
\numberwithin{figure}{section}
\numberwithin{table}{section}

We continue our study \cite{firo23} of the global dynamics for the scalar reaction-advection-diffusion equation
	\begin{equation}
	u_t = u_{xx} + f(x,u,u_x)\,.
	\label{PDE}
	\end{equation}
Subscripts $t,x$ indicate partial derivatives.
To be specific, we consider solutions $u=u(t,x) \in \mathbb{R}$ on the unit interval $0<x<1$, with Neumann conditions $u_x=0$ at the boundaries $x=0,1$.
Equilibria of \eqref{PDE}, i.e.~time-independent solutions $u(t,x)=v(x)$, equivalently satisfy the ``pendulum'' equation  
	\begin{equation}
	0 = v_{xx} + f(x,v,v_x)\,,
	\label{ODE}
	\end{equation} 
albeit as a Neumann boundary value problem in the spatial variable $x$.

The mathematical literature on reaction-diffusion equations alone, as refereed in Zentralblatt under MSC 35K57, has grown to more than 15,000 entries \cite{zb}.
We have already provided detailed mathematical background, and a survey of applications; see \cite{firo23} and the many references there.
We only recall some of the most important facts, for the convenience of our readers.
Our main goal is to proceed from basic information on the ODE boundary value problem \eqref{ODE}, and derive detailed structural results on the global attractor of the PDE \eqref{PDE}.
However, we do not start from given nonlinearities $f$, as the primary object.
Due to the (spatial) chaoticity of \eqref{ODE} on large $x$-intervals, even the task to determine all equilibria can in fact be prohibitively messy.
Quite surprisingly, therefore, it is possible to characterize the class of all ODE equilibrium ``configurations'', qualitatively, by certain permutations $\sigma$, as introduced by Fusco and Rocha \cite{furo91}.
See \eqref{hdef} and \eqref{permdef} below.
Suffice it to assert here that any of those permutations $\sigma$ can be represented by an open class of dissipative nonlinearities $f$.
The 3-nose meanders, which we continue to study in the present paper, are then defined in terms of Fusco's permutations $\sigma$, rather than by their associated (classes of) nonlinearities $f$.

We assume the PDE solution semigroup $u(t,\cdot)$ generated by the nonlinearity $f\in C^2$ to be \emph{dissipative}: any solution $u(t,\cdot)$ exists globally in forward time $t\geq 0$, and eventually enters a fixed large ball in a suitable Sobolev space $X$ which contains $C^1 ([0,1], \mathbb{R})$.
(The sign condition $f(x,u,0)\cdot u<0$, for large $|u|$, together with subquadratic growth of $f(x,u,p)$ in $|p|$, for example, is sufficient for that.)
For large times $t\rightarrow\infty$, any large ball in $X$ then limits onto the same maximal compact and invariant subset $\mathcal{A}=\mathcal{A}_f$ of $X$ which is called the \emph{global attractor}. 
In general, the global attractor $\mathcal{A}$ consists of all solutions $u(t,\cdot)$ which exist globally, for all positive and negative times $t\in\mathbb{R}$, and which remain uniformly bounded in $X$.
In general, $\mathcal{A}$ therefore contains any equilibria, heteroclinic orbits, basin boundaries, or more complicated recurrence which might arise.
Let $\mathcal{E} \subseteq \mathcal{A}$ denote the set of equilibria.
We assume all equilibria are \emph{hyperbolic}.

\begin{figure}[p!]
\centering \includegraphics[width=0.9\textwidth]{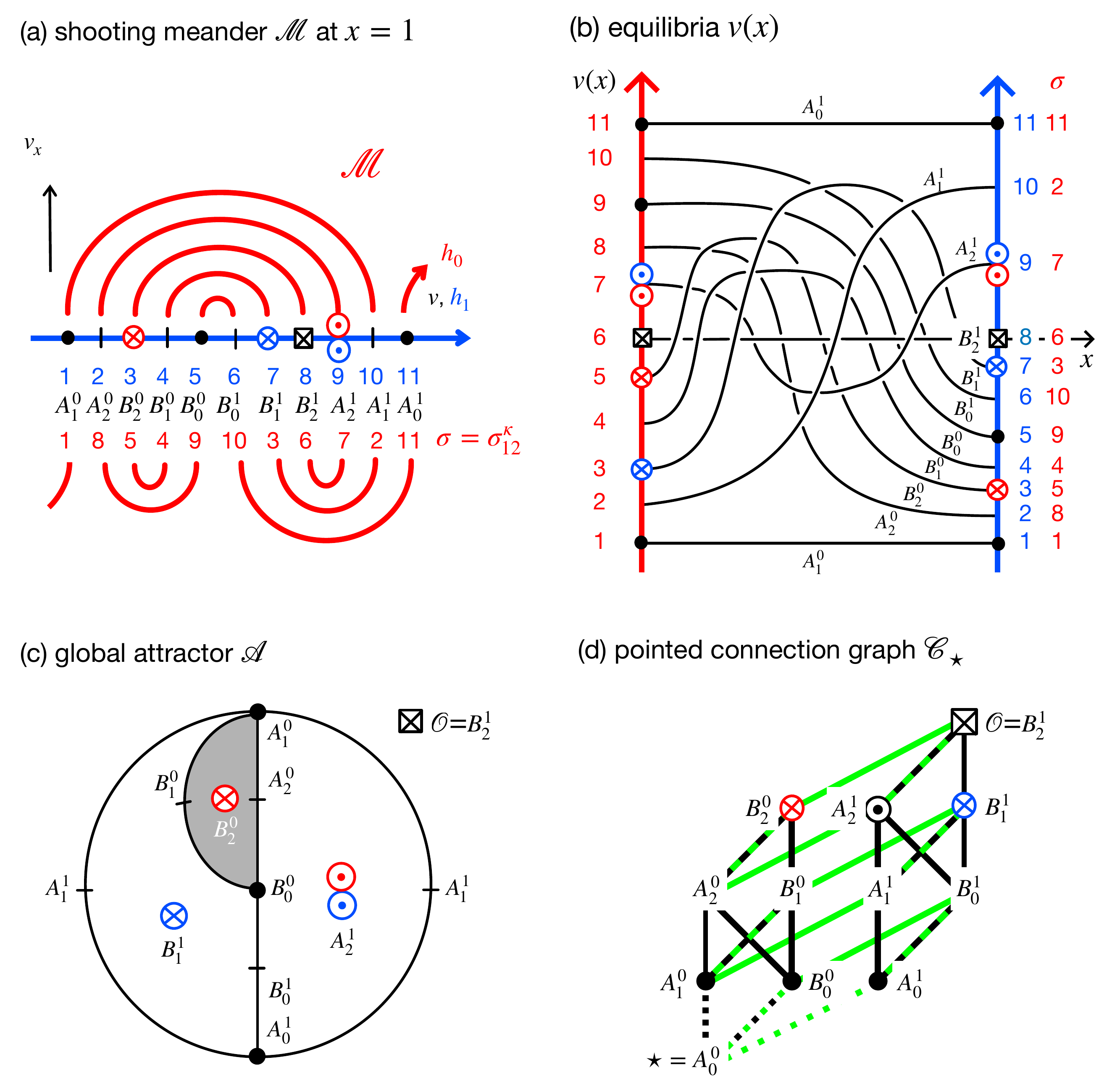}
\caption{\emph{
Example of a Sturm 3-ball global attractor $\mathcal{A}= \mathrm{clos}\ W^u(\mathcal{O})$. 
Equilibria \hbox{are labeled} as $\mathcal{E}=
\{A^0_1,A^0_2,B^0_2,B^0_1,B^0_0,B^1_0,B^1_1,B^1_2=\mathcal{O},A^1_2,A^1_1,A^1_0\}$, along the $h_1$-axis, i.e.~ordered by their boundary values $v(x)$ at $x=1$. 
Note the Morse indices $i(L^j_k)=j+k$, for tag $L=B$, and $i=j+k-1$ for tag $L=A$.
Black dots mark sinks $i=0$, and small annotated circles (red, blue) indicate $i=2$.
The previous papers \cite{firo3d-1, firo3d-2} established the equivalence of the viewpoints (a)--(c).\\
(a) The Sturm meander $\mathcal{M}$ of the global attractor $\mathcal{A}$. 
The meander $\mathcal{M}$ is a stylized representation of the curve $a \mapsto (v,v_x)$, at $x=1$, which results from Neumann initial conditions $(v,v_x)=(a,0)$, at $x=0$, by shooting via the equilibrium ODE \eqref{ODE}. 
Intersections of the meander with the horizontal $v$-axis indicate equilibria $v(x)$. 
Transversality of intersections is equivalent to hyperbolicity of $v$.\\
(b) Spatial profiles $x\mapsto v(x)$ of the equilibria $v \in \mathcal{E}$. 
Note the different orderings of $v(x)$, by $h_0$ at the left boundary $x=0$ (red), but by $h_1$ (blue) and by the Sturm permutation $\sigma = h_0^{-1}h_1$ (red) at the right boundary $x=1$. 
The same orderings characterize the meander $\mathcal{M}$ in (a).\\
(c) The Thom-Smale or Sturm complex $\mathcal{S}$ of the boundary 2-sphere $\Sigma^2=\partial \mathcal{A}=\partial W^u(\mathcal{O})$.
The right and left boundaries $A^0_1A^1_1A^1_0$ both indicate the unstable manifold of $A^1_1$ and have to be identified with each other. 
Note the 2-cells given by the unstable manifolds of the three equilibria $B^1_1,B^0_2,A^1_2$ of Morse index $i=2$.
The Chafee-Infante 2-cell of $B^0_2$ is shaded gray.\\
(d) The connection graph $\mathcal{C}$ of $\mathcal{A}$.
Vertices are the equilibria, ranked by rows of equal Morse index.
Edges (green or black) are directed downwards and indicate heteroclinic orbits $\leadsto$ between Morse-adjacent equilibria.
Green rhombics indicate three slanted Chafee-Infante stacks of height 2 each;
black edges mark two vertical Chafee-Infante stacks of height 3; 
dashed black/green edges are shared between different stacks.
See sections \ref{ChIn} and \ref{Conn} for further details.
}}
\label{3ball}
\end{figure}

Two additional structures help to describe $\mathcal{A}$, in our specific setting \eqref{PDE}.
First, \eqref{PDE} possesses a \emph{Lyapunov~function}, alias a variational or gradient-like structure, under separated boundary conditions;  see \cite{ze68, ma78, mana97, hu11, fietal14, lafi18, labe22}. 
Therefore, the time-invariant global attractor consists of equilibria and of solutions $u(t, \cdot )$, $t \in \mathbb{R}$, with forward and backward limits, i.e.
	\begin{equation}
	\underset{t \rightarrow -\infty}{\mathrm{lim}} u(t, \cdot ) = v_1\,,
	\qquad
	\underset{t \rightarrow +\infty}{\mathrm{lim}} u(t, \cdot ) = v_2\,.
	\label{het}
	\end{equation}
In other words, the $\alpha$- and $\omega$-limit sets of $u(t,\cdot )$ are two distinct hyperbolic equilibria $v_1$ and $v_2$.
We call $u(t, \cdot )$ a \emph{heteroclinic} or \emph{connecting} orbit, or \emph{instanton},  and write $v_1 \leadsto v_2$ for such heteroclinically connected equilibria. 
See fig.~\ref{3ball}(c),(d) for a modest 3-ball example with $N=11$ equilibria.
Although the gradient-like structure persists for other separated boundary conditions, the possibility of rotating waves shows that it may fail under periodic boundary conditions.

The second structure is a \emph{Sturm nodal property}, which we express by the \emph{zero number} $z$.
This nodal property justifies to call the global attractors $\mathcal{A}$ \emph{Sturm}.
Let $0 \leq z (\varphi) \leq \infty$ count the number of strict sign changes of continuous spatial profiles $\varphi : [0,1] \rightarrow \mathbb{R}, \, \varphi \not\equiv 0$.
For any two distinct solutions $u^1$, $u^2$ of \eqref{PDE}, the zero number
	\begin{equation}
	t \quad \longmapsto \quad z(u^1(t, \cdot ) - u^2(t, \cdot ))\ \searrow
	\label{zdrop}
	\end{equation}
is then nonincreasing with time $t$, for $t\geq0$, and finite for $t>0$.
Moreover $z$ drops strictly, with increasing $t>0$, at any multiple zero of the spatial profile $x \mapsto u^1(t_0 ,x) - u^2(t_0 ,x)$; see \cite{an88}.
This remains true under other separated or periodic boundary conditions.
See Sturm \cite{st1836} for the linear autonomous variant, and \cite{ma82} for the PDE revival.. 

The consequences of the Sturm nodal property \eqref{zdrop} for the nonlinear dynamics of \eqref{PDE} are enormous.
For example, Morse-Smale transversality, and hence structural stability, hold automatically, just given hyperbolicity of equilibria \cite{he85, an86}.
Already Sturm observed that all eigenvalues $\mu_0>\mu_1>\ldots$ of the PDE linearization of \eqref{PDE} at any equilibrium $v$ are real and algebraically simple.
In fact $z(\varphi_j)=j$, for the eigenfunction $\varphi_j$ of $\mu_j$.
The \emph{Morse index} $i(v)$ of $v$ then counts the number of unstable eigenvalues $\mu_j>0$.
The Morse index $i(v)$ is the dimension of the unstable manifold $W^u(v)$ of $v$.

In a series of papers, based on the zero number, we have given a purely combinatorial description of Sturm global attractors $\mathcal{A}$; see \cite{firo96, firo99, firo00}.
Define the two \emph{boundary orders} $h_0, h_1$: $\lbrace 1, \ldots, N \rbrace \rightarrow \mathcal{E}$ of the equilibria such that
	\begin{equation}
	h_\iota (1) < h_\iota (2) < \ldots < h_\iota (N) \qquad \mathrm{at}
	\qquad x=\iota \in \{0,1\}\,.
	\label{hdef}
	\end{equation}
See fig.~\ref{3ball}(a),(b) for an example with $N=11$ equilibrium profiles, $\mathcal{E} = \{1,\ldots,11\}, \ h_0 = \mathrm{id},\ h_1 = (1\ 8\ 5\ 4\ 9\ 10\ 3\ 6\ 7\ 2\ 11)$.	
The general combinatorial description of Sturm global attractors $\mathcal{A}$ is based on Fusco's groundbreaking \emph{Sturm~permutation} $\sigma \in S_N$\,, defined in \cite{furo91} as
	\begin{equation}
	\sigma:= h_0^{-1} \circ h_1\,.
	\label{permdef}
	\end{equation}
Note how $\sigma$ is trivially invariant under any bijective relabeling $\Lambda: \mathcal{E}_1 \rightarrow \mathcal{E}_2$ of equilibria, because $(\Lambda h_0)^{-1}(\Lambda h_1)=h^{-1}_0 h_1$. 
Much less trivially, Sturm attractors of dissipative nonlinearities with the same Sturm permutation $\sigma$ are $C^0$ orbit-equivalent \cite{firo00}.

Already in \cite{furo91}, the following explicit recursions have been derived for the Morse indices $i_j:=i(h_0(j))$ along the meander:
	\begin{equation}
	\begin{aligned}
	i_1 &=  i_N = 0\,;\\
	i_{j+1} &= i_j+(-1)^{j+1}\,
	\mathrm{sign}\, (\sigma^{-1}(j+1)-\sigma^{-1}(j))\,.\\
	\end{aligned}
	\label{i}
	\end{equation}
The zero numbers, $z_{jk} := z(h_0(j)-h_0(k))\geq 0$ for $j\neq k$, are given recursively by
	\begin{equation}
	\begin{aligned}
	z_{kk} &:= i_k\,;\qquad\qquad\qquad \\
	z_{1k} &\phantom{:}= z_{Nk}=0\,;\\
	z_{k+1,k} &\phantom{:}= \min\{i_k\,,i_{k+1}\}\,;\\
	z_{j+1,k} &\phantom{:}= \scalebox{0.97}[1.0]
	{$ z_{jk} + \tfrac{1}{2}(-1)^{j+1}
	\cdot\big[ \mathrm{sign}\,\big(\sigma^{-1}(j+1)-\sigma^{-1}(k)\big)-\mathrm{sign}\,\big(\sigma^{-1}(j)-\sigma^{-1}(k)\big)\big].
	$}
	\end{aligned}
	\label{z}
	\end{equation}


Using a shooting approach to the ODE boundary value problem \eqref{ODE}, the Sturm~permutations $\sigma \in S_N$ have been characterized, purely combinatorially, as \emph{dissipative Morse meanders} in \cite{firo99}.
Here the \emph{dissipativeness} property, abstractly, requires fixed $\sigma(1)=1$ and $\sigma(N)=N$.
In fact, the shooting meander emanates upwards, towards $v_x>0$, from the leftmost (or lowest) equilibrium at $\sigma(1)=1$, and terminates from below, $v_x<0$, at x=1.
The \emph{meander} property requires the formal path $\mathcal{M}$ of alternating upper and lower half-circle arcs defined by the permutation $\sigma$, as in fig.~\ref{3ball}(a), to be Jordan, i.e.~non-selfintersecting.
For dissipative meanders, the recursion \eqref{i}, and $i_1=0$, define all Morse numbers $i_j$\,.
Note how $j$ and $i_j$ are always of opposite parity, $\mathrm{mod}\ 2$.
In particular, $N$ is odd, and $i_N$=0 follows automatically.
The \emph{Morse} property, finally, requires nonnegative Morse indices $i_j\geq0$ in \eqref{i}, for all $j$.
For brevity, we also use the term \emph{Sturm meanders}, for dissipative Morse meanders.

For a simple recipe to determine the Morse property of a meander, \emph{the Morse number increases by 1, along any right turning meander arc, but decreases by 1 along left turns.}
This holds, independently, for upper and lower meander arcs, and remains valid even when the proper orientation of the arc is reversed; see \eqref{i}.
For examples see figs.~\ref{3ball}, \ref{ci3}.


In \cite{firo96} we have shown how to determine which equilibria $v_1,v_2$ possess a heteroclinic orbit connection \eqref{het}, explicitly and purely combinatorially from dissipative Morse meanders $\sigma$. 
In the elegant formulation of Wolfrum  \cite{wo02},
\begin{equation}\label{wolfrum}
v_1\leadsto v_2 \quad\Longleftrightarrow\quad  v_1, v_2\  \textrm{are}\ z\textrm{-adjacent, and}\ i(v_1)>i(v_2)\,;
\end{equation}
see also the comment in the appendix of \cite{firo3d-2}.
Here equilibria $v_1\neq v_2$ are called $z$-\emph{adjacent}, if there does \emph{not} exist any \emph{blocking equilibrium} $w$ strictly between $v_1$ and $v_2$, at $x=0$ (or, equivalently, at $x=1$), such that
\begin{equation}\label{block}
z(v_1-w)=z(w-v_2)=z(v_1-v_2).
\end{equation}
With \eqref{wolfrum}, all heteroclinic orbits then follow from \eqref{i} and \eqref{z} above.

Clearly, any heteroclinic orbit $u(t,.): v_1\leadsto v_2$ implies adjacency: by \eqref{zdrop}, any blocking $w$ would force $z(u(t,.)-w)$ to drop strictly at the Neumann boundary $x=0$, for some $t=t_0$.
This contradicts the equal values \eqref{block} of $z$ at the limiting equilibria $v_1,v_2$ of $u$, for $t\rightarrow\pm\infty$.

As a trivial corollary, for example, we conclude $v_1\leadsto v_2$, for neighbors $v_1,v_2$ on any boundary order $h_\iota$\,. Here we label $v_1,v_2$ such that $i(v_1)=i(v_2)+1$; see \eqref{i}.
For in-depth analysis and many more examples see \cite{furo91, firo14, firo15, firo3d-3, rofi21}. 

We summarize the above heteroclinic structure in the directed graded \emph{connection graph} $\mathcal{C}$.
See fig.~\ref{3ball}(d) for an example.
The connection graph is graded by increasing levels of the Morse indices $i$ of its equilibrium vertices, as a ranking function.
Directed edges are the heteroclinic orbits $v_1\leadsto v_2$ running downwards between certain equilibria of adjacent Morse index.
Uniqueness of such heteroclinic orbits, given $v_1,v_2$, had already been observed in lemma 3.5 of \cite{brfi89}.
The connection graph in fact encodes all heteroclinic orbits, due to transitivity of the relation $\leadsto$ (by Morse-Smale transversality), and a cascading principle \cite{brfi89, fi94, firo96}.

For convenience, we also define an augmentation $\mathcal{C}_\star$ of $\mathcal{C}$, which we call the \emph{pointed connection graph}, or the \emph{connection graph with a distinguished vertex} $\star$\,.
This is reminiscent of Conley theory, where the (empty!) exit set of an attractor is usually collapsed into a ``distinguished'' point \cite{con, mis}.
We simply add a ``distinguished'' vertex $\star$ at ``Morse'' level $i=-1$, with an edge towards $\star$ from every sink equilibrium $v$, i.e.~from every vertex $v$ with Morse index $i(v)=0$.

More geometrically, the global attractor $\mathcal{A}$ is a signed regular cell complex $\mathcal{S}$.
We call $\mathcal{S}$ the \emph{Sturm complex} or, because the cells are given by the unstable manifolds of the hyperbolic equilibria, the \emph{Thom-Smale complex}; see \cite{firo20, rofi21}.
We call $d=\dim \mathcal{A} := \max_{v \in \mathcal{E}} i(v)$ the \emph{dimension} of $\mathcal{A}$ or of the cell complex $\mathcal{S}$.
Then at least one equilibrium $\mathcal{O}$ has maximal Morse index $i(\mathcal{O})=d$, i.e. $i(v) \leq d$ for all other Morse indices.
If $\mathcal{A} = \mathrm{clos}\,W^u(\mathcal{O})$ is the closure of a single $d$-cell, then the Thom-Smale complex turns out to be a closed $d$-ball \cite{firo15}.
We call this case a \emph{Sturm $d$-ball}.
See fig.~\ref{3ball}(c) for the Sturm 3-ball associated to the meander in fig.~\ref{3ball}(a).

In the present paper we discuss Sturm attractors which arise from Sturm meanders with at most three \emph{noses}.
Noses are defined by $j\in\{1,\ldots,N-1\}$ such that
\begin{equation}
\label{nose}
\sigma(j+1)=\sigma(j)\pm 1.
\end{equation}
For a 3-nose meander see fig.~\ref{3ball}(a) again.
The simplest case, of just two noses, is called the \emph{Chafee-Infante attractor}.
In 1974, this case arose for cubic nonlinearities \eqref{cubic}; see \cite{chin74}.
Chafee-Infante components turn out to be fundamental building blocks of 3-nose connection graphs.
We therefore revisit and review this case in section \ref{ChIn}.

We now present our main results on the general case of \emph{primitive 3-nose meanders} $\mathcal{M}_{pq}$\,.
These meanders are characterized by $p$ nested arcs around their upper \emph{left} innermost nose arc,
above the horizontal axis, and by $q$ nested arcs around their upper \emph{right} innermost nose arc.
Below the horizontal axis, the only remaining nose is centered as the innermost of the remaining $p+q$ lower arcs.
Since all lower arcs are nested, we also call that configuration a (lower) \emph{rainbow}.
See fig.~\ref{qnest}(b) for a general template.

In theorem 2.1 of \cite{firo23}, we have shown that $\mathcal{M}_{pq}$ is a dissipative meander if, and only if, $(p-1,q+1)$ are coprime and $p\geq2$.
Moreover, such meanders fail to be Morse, for $p\neq r(q+1)$.
Therefore, we focus on the 3-nose cases
\begin{equation}
\label{prq}
p=r(q+1)
\end{equation} 
Note that $p-1$ and $q+1$ are then coprime, indeed, because $r(q+1)-(p-1)=1$.

The trivial case $r=q=1$ has been illustrated in fig.~3.2(a)-(c) of \cite{firo23}.
We therefore assume $r,q\geq 1, \ rq>1$, for the rest of this paper.
Only for $r=1$ have the next four theorems been proved in \cite{firo23}.
The corollaries have been proved, provided the theorems hold.

Our first theorem shows that, conversely, all cases \eqref{prq} lead to Morse meanders $\mathcal{M}_{pq}$.

\begin{thm}\label{krmuthm}
For $p=r(q+1)$, let $\mu_{rq}(i)$ count the vertices with Morse number $i$, in the dissipative meander $\mathcal{M}_{pq}$\,. 
Then, for $r,q\geq 1, \ rq>1$, the nonzero Morse counts are given by 
\begin{equation}\label{krmu}
\mu_{rq}(i)\ =\ 
\begin{cases}
      3+2i,& \textrm{for}\quad 0\leq i< \min\{r,q\};\\
      2+2\min\{r,q\},& \textrm{for}\quad \min\{r,q\}\leq i < \max\{r,q\};\\
      2(r+q)+1-2i,& \textrm{for}\quad \max\{r,q\}\leq i\leq r+q.
\end{cases}
\end{equation}
In particular, all such meanders $\mathcal{M}_{pq}$ are Sturm.
\end{thm}

\begin{cor}\label{krmucor}
The Morse count functions $i\mapsto \mu_{rq}(i)$ have the following symmetry properties:
\begin{enumerate}[(i)]
  \item  up to ordering, the subscript set $\{r,q\}$ is determined by $\mu_{rq}$\,;
  \item conversely, the subscript set determines $\mu_{rq}=\mu_{qr}$\,;
  \item for all $0\leq i < r+q$, we have $\mu_{rq}(i) = \mu_{rq}(r+q-1-i)$.
\end{enumerate}
\end{cor}

\begin{defi}\label{Akrdefi}
For the Sturm entourage of primitive 3-nose meanders $\mathcal{M}_{r(q+1),q}$\,, we denote the associated \emph{primitive
Sturm permutation} as $\sigma_{rq}$\,, the \emph{primitive Sturm attractor} as $\mathcal{A}_{rq}$\,, and the \emph{primitive connection graph} as $\mathcal{C}_{rq}$\,.
\end{defi}

For our next result we have to recall the notion of \emph{trivial equivalences} from section 3 of \cite{firo23}.
In the PDE context \eqref{PDE}, these consist of the Klein 4-group $\langle \kappa,\varrho \rangle$ with commuting involutive generators
	\begin{align}
	(\kappa u)(x) &:= -u(x)\,; \label{rot}\\
	(\varrho u)(x) &:= u(1-x)\,.\label{inv}
	\end{align}
The involution \eqref{rot} therefore simply rotates $\mathcal{M}\subset \mathbb{R}^2$ by $180^\circ$, i.e. 
\begin{equation}
\label{rotM}
\mathcal{M}^\kappa:=-\mathcal{M}\,.
\end{equation}
The orientation of the meander curve, however, is reversed.
Abusing notation slightly, let $\kappa$ also denote the flip permutation
	\begin{equation}\label{flip}
	\kappa (j) := N+1-j
	\end{equation}
on $j \in \{ 1, \ldots, N\}$.
Then 
the meander rotation \eqref{rotM} leads to conjugation 
	\begin{equation}\label{rotsig}
	\sigma^\kappa = \kappa \sigma \kappa
	\end{equation}	
by the flip involution \eqref{flip}.

Spatial reversal $\varrho$ of $x$, in contrast, interchanges the boundaries $x=\iota \in \{0,1\}$.
The effect on Sturm permutations $\sigma = h_0^{-1}\circ h_1$ is inversion:
	\begin{equation}\label{invsig}
	\sigma^\varrho = \sigma^{-1}\,.
	\end{equation}

\begin{thm}\label{rqqrthm}
The primitive Sturm permutations $\sigma_{rq}$ and $\sigma_{qr}$ are trivially equivalent under the involutive product $\kappa\varrho$ of \eqref{rotsig} and \eqref{invsig}.
In symbols,
\begin{equation}\label{rqqr}
\sigma_{qr}=\sigma^{\kappa\varrho}_{rq}=\kappa(\sigma_{rq})^{-1}\kappa.
\end{equation} 
\end{thm}

\begin{cor}\label{rqqrCor}
The primitive 3-nose Sturm attractors $\mathcal{A}_{rq}$ and $\mathcal{A}_{r'q'}$ are orbit equivalent if, and only if, their subscript sets coincide, up to ordering.
In fact $\mathcal{A}_{rq}$ and $\mathcal{A}_{qr}$ are trivially equivalent, under the involutive product $\kappa\varrho$ of \eqref{rot} and \eqref{inv}.
\end{cor}

By a simple rotation $\kappa$ of the meander $\mathcal{M}$, for example, the 3-ball attractor of fig.~\ref{3ball} is trivially equivalent to the simple case $r=1,\ q=2$ with Sturm permutation $\sigma_{12}$.

Note the Morse count $\mu_{rq}(r+q)=1$ at maximal $i=r+q$ in \eqref{krmu}. 
Let $\mathcal{O}$ denote that unique equilibrium in $\mathcal{A}_{rq}$ of maximal Morse index $i(\mathcal{O})=r+q=\dim \mathcal{A}_{rq}$\,.

\begin{thm}\label{krballthm}
The primitive Sturm attractor $\mathcal{A}_{rq}$ is the closure of the unstable manifold of the single equilibrium $\mathcal{O}\in \mathcal{A}_{rq}$\,.
I.e., $\mathcal{A}_{rq}$ is a Sturm ball of dimension $r+q$.
\end{thm}

Let $\Sigma^{r+q-1} := \partial\mathcal{A} = \partial W^u(\mathcal{O}) := \mathrm{clos} W^u(\mathcal{O}) \setminus W^u(\mathcal{O})$ denote the invariant boundary $(r+q-1)$-sphere of the $(r+q)$-dimensional Sturm ball $\mathcal{A}_{rq}$\,.
Quite surprisingly, the connection graph $\mathcal{C}_{rq}$\,, restricted to $\Sigma^{r+q-1}$, then turns out to be time-reversible; see theorem \ref{krrevthm} below.
Although this is also true in the Chafee-Infante case, it is a quite unexpected phenomenon in general parabolic diffusion equations -- which most of us would rightly consider \emph{the} paradigm of irreversibility.

\emph{Time reversibility} in its strongest form means the existence of an involutive \emph{reversor} $\mathcal{R}: \Sigma \rightarrow \Sigma$ which reverses the time direction of PDE orbits of \eqref{PDE}, on a ``large'' invariant subset $\Sigma\subset\mathcal{A}$.
Restricted to equilibria $v_1,v_2\in \mathcal{E} \subset \Sigma$, strong reversibility implies the weaker statement
\begin{equation}
\label{rev}
v_1\leadsto v_2 \quad\Longleftrightarrow\quad \mathcal{R}v_2 \leadsto \mathcal{R}v_1
\end{equation}
on $\Sigma$. 
In other words, the reversor $\mathcal{R}$ induces an automorphism of the connection di-graph $\mathcal{C}|_\Sigma$\,, which reverses heteroclinic edge orientation.
We then call the connection graph \emph{time-reversible}.

\begin{thm}\label{krrevthm}
The pointed connection graphs $\mathcal{C}_{rq\star}$ are time-reversible.
\end{thm}

Consider the pointed connection graph $\mathcal{C}_\star = \mathcal{C}^\kappa_{12\star}$ of fig.~\ref{3ball}(d), for example.
Then edge-reversing invariance of $\mathcal{C}_\star$ under rotation by $180^\circ$ illustrates reversibility \eqref{rev} under the explicit reversor
\begin{equation}
\label{3ballrev}
\mathcal{R}: \quad A^k_j \longleftrightarrow B^{1-k}_{2-j}
\end{equation}
of the equilibria on the 2-sphere $\Sigma:=\partial W^u(\mathcal{O})$ of fig.~\ref{3ball}(c).

As in the general case, the artificial edges towards the artificial vertex $\star$ at ``Morse'' index $i=-1$ do not signify heteroclinic orbits. 
Since the involutive reversor $\mathcal{R}$ on $\mathcal{C}_{rq}$ swaps vertices of Morse levels $i$ and $r+q-1-i$, we have $\mathcal{O} = \mathcal{R}\,\star$ as the reversor counterpart of $\star$, in general.
The actual heteroclinic orbits emanating from $\mathcal{O}$, however, cannot be matched by actual counterpart PDE orbits towards $\star$, under the reversor $\mathcal{R}$.
Therefore reversibility \eqref{rev} for actual heteroclinic orbits can only be asserted among equilibria at Morse levels $i=0,\ldots i(\mathcal{O})-1$, i.e.~on the flow-invariant boundary sphere $\Sigma^{r+q-1}=\partial\mathcal{A}_{rq}=\partial W^u(\mathcal{O})$ of the Sturm ball $\mathcal{A}_{rq} = \mathrm{clos}\ W^u(\mathcal{O})$.
That reversibility on the boundary sphere $\Sigma^{r+q-1}$, of course, is a much deeper reason for the symmetry of the Morse count function $i\mapsto \mu_{rq}(i)$, for $0\leq i < r+q$, which we have already encountered in corollary \ref{krmucor}(iii).

The remaining paper is organized as follows.
We briefly recall the Chafee-Infante paradigm, in section \ref{ChIn}, to settle notation.
Section \ref{Sigma} presents explicit expressions for the meander permutations $\sigma_{rq}$ of the dissipative meanders $\mathcal{M}_{pq}$ with $p=r(q+1)$.
This provides a proof of claim \eqref{rqqr}, i.e.~of the trivial equivalence of $\sigma_{qr}$ and $\sigma_{rq}$ stated in theorem \ref{rqqrthm}.
We also label the vertices of $\mathcal{M}_{pq}$ for a convenient description of the pointed connection graph $\mathcal{C}_{rq\star}$ in terms of pointed Chafee-Infante stacks $\mathcal{C}_{r\star}$ and $\mathcal{C}_{q\star}$\,.
See \eqref{ciequi}, \eqref{cis}.
The connection graph itself will be described as a lattice of Chafee-Infante stacks, in our main theorem \ref{connthm}.
Sections \ref{krmupf}--\ref{krrevpf} will then show how all remaining theorems \ref{krmuthm}, \ref{rqqrthm}, \ref{krballthm}, and \ref{krrevthm} follow from the main theorem \ref{connthm}.
Theorem \ref{connthm}, in turn, is proved in section \ref{Rec} by a double recursion on, both, $r$ and $q$.
In section \ref{Dis}, we embark on  a discussion of non-Morse variants of our results, with a focus on time reversibility of 3-nose connection graphs involving suspensions.

\textbf{Acknowledgment.}

We dedicate this paper to Giorgio Fusco, a most amiable Gentiluomo, for his long friendship and his many generously shared refined insights, mathematical and otherwise. 

Support by FCT/Portugal under UID/MAT/04459/2019 and UIDB/04459/2020 is gratefully acknowledged.
This work has also been most generously supported by the Deutsche Forschungsgemeinschaft, Collaborative Research Center 910 \emph{``Control of self-organizing nonlinear systems: Theoretical methods and concepts of application''} under project A4: \emph{``Spatio-temporal patterns: control, delays, and design.''}

\section{Chafee-Infante stacks}\label{ChIn}

In this section we recall some facts on the sequence $\mathcal{M}_d\,,\ d\geq1$, of Sturm meanders with two noses and $2d$ arcs.
This class was  studied by Chafee and Infante \cite{chin74} in the guise of PDE \eqref{PDE} with symmetric cubic nonlinearity $f$ and parameter $\lambda$,
\begin{equation}\label{cubic}
f=\lambda^2 u (1-u^2), \qquad(d-1)\pi<\lambda<d\pi.
\end{equation}
The 2-nose Sturm meander $\mathcal{M}_d$ then arises by a time map analysis of the Duffing type ODE \eqref{ODE}.
We therefore attach the names \emph{Chafee-Infante} to the meander $\mathcal{M}_d$ and its entourage of Sturm permutation $\sigma_d$\,, Sturm attractor $\mathcal{A}_d$\,, connection graph $\mathcal{C}_d$\,, and pointed connection graph $\mathcal{C}_{d\star}$\,.
We call the pointed version $\mathcal{C}_{d\star}$ with $2(d+1)$ vertices a \emph{Chafee-Infante stack} (CIS) of \emph{height} $d+1$.
For further details on the following facts and remarks we refer to section 4 of our prequel \cite{firo23}, as well as to \cite{he85, fi94, firo20}.

The Chafee-Infante meander $\mathcal{M}_d$ is defined by $d$ nested upper arcs, and $d$ nested lower arcs, correspondingly shifted by one vertex to the right.
In other words, $\mathcal{M}_d$ consists of an upper and a lower rainbow, each with $d$ nested arcs.
Explicitly, the associated Sturm permutation is
\begin{equation}\label{cisigma}
\sigma_d(j)  \ =\  
\begin{cases}
      \  j, &\textrm{ for odd } j,  \\
      \ N+1-j, &\textrm{ for even } j.
\end{cases} 
\end{equation}
Here $j=1,\ldots,N$ enumerates the $N=2d+1$ equilibria.
Note invariance $\sigma_d^\kappa = \sigma_d^\varrho = \sigma_d$ under the trivial equivalences \eqref{rotsig} and \eqref{invsig}.

\begin{wrapfigure}[21]{l}{0.5\textwidth}
\centering 
\includegraphics[width=\linewidth]{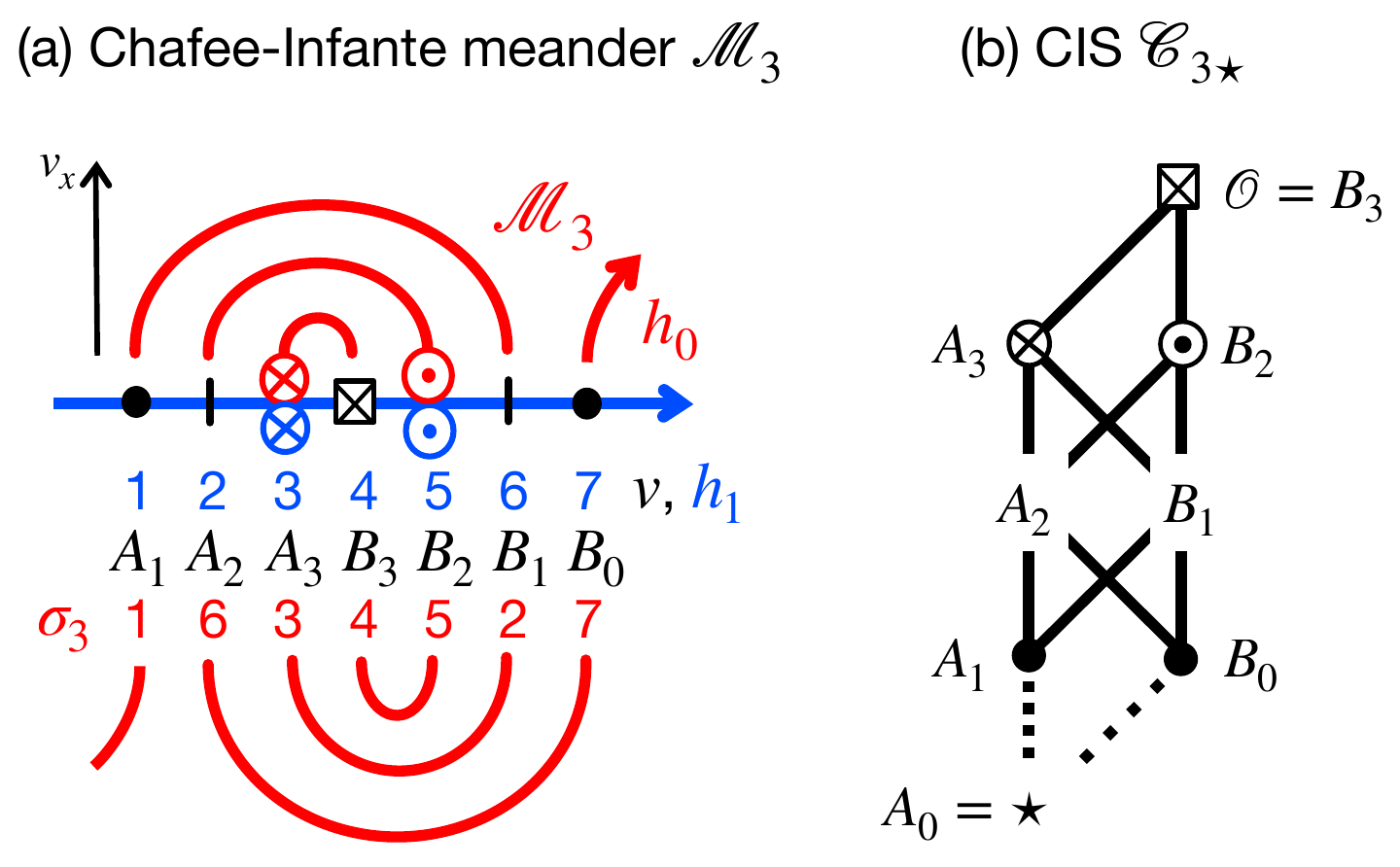}
\caption{\emph{
The Chafee-Infante attractor of dimension $d=3$; 
see PDE \eqref{PDE} with cubic nonlinearity \eqref{cubic}.\\
(a) The associated stylized Sturm meander $\mathcal{M}_3$. See \eqref{cisigma}, \eqref{ciequi}.\\
(b) The pointed connection di-graph $\mathrm{CIS}\ \mathcal{C}_{3\star}$\,, ranked by Morse index; see \eqref{cis}. 
Time reversibility in the boundary sphere of $\mathcal{C}_3$ appears as an automorphism of $\mathcal{C}_{3\star}$ under the $180^\circ$ rotation by the reversor $\mathcal{R}: A_j \leftrightarrow B_{3-j}\,$; see  \eqref{cirev}.
Note how the involution $\mathcal{R}$ reverses the Morse ranking, and hence the downward edge orientations.
}}
\label{ci3}
\end{wrapfigure}
To describe the pointed CIS $\mathcal{C}_{d\star}$ of height $d+1$, we label the $N=2d+1$ equilibrium vertices along the horizontal axis as
\begin{equation}\label{ciequi}
h_1(j) \ =:\ 
\begin{cases}
      \ A_{j}\,, & \textrm{for } \qquad 1\leq j \leq d, \\
      \ B_{N-j}\,, & \textrm{for } \,d+1 \leq j \leq N.
\end{cases}
\end{equation}
We also adjoin the distinguished vertex $A_0:= \star$\,.
The CIS is then given by
\begin{equation}\label{cis}
A_{j+1}\,,\, B_j \quad \leadsto \quad A_j\,,\, B_{j-1} \,.
\end{equation}
This means that each vertex on the left of \eqref{cis} possesses a directed edge to each vertex on the right, for all appropriate subscripts. 
Subscripts range from $1$ to $d$, on the left, and from $0$ to $d-1$ on the right.
See fig.~\ref{ci3} for an illustration of the Chafee-Infante meander $\mathcal{M}_d\,,\ d=3$, and the associated CIS $\mathcal{C}_{d\star}$ of height $d+1$ with the $2(d+1)$ vertices $A_0,\ldots,A_d$ and $B_0,\ldots,B_d$.

In particular, the Morse numbers of any vertex $L_j$, with tags $L = A,B$, are
\begin{equation}
\label{iLj}
i(L_j)\ =\ \begin{cases}
      \ j-1& \text{for } L=A, \\
      \ j& \text{for } L=B.
\end{cases}
\end{equation}
Therefore $\mathcal{O}:=B_d$ is the unique equilibrium of maximal Morse index $d$.
Moreover, $B_d$ connects to any other vertex $E$ by a di-path in the CIS, i.e.~the equilibrium $B_d \leadsto E$ connects to any other equilibrium $E$, heteroclinically.
In particular, $\mathcal{A}_d= \mathrm{clos}\, W^u(\mathcal{O})$ is a Sturm $d$-ball and $\partial W^u(\mathcal{O}) = \Sigma^{d-1}$ is a sphere of dimension $d-1$.

This also proves the following equivalent extremal characterizations of the Chafee-Infante attractor $\mathcal{A}_d$, among all Sturm attractors \cite{fi94,firo20}\,:
\begin{description}
  \item[max:] \emph{Among all Sturm attractors with $N=2d+1\geq 3$ equilibria, $\mathcal{A}_d$ is the unique Sturm attractor with the maximal possible dimension $d$;}
  \item[min:] \emph{Among all Sturm attractors of dimension $d\geq 1$, $\mathcal{A}_d$ is the unique Sturm attractor with the smallest possible number $N=2d+1$ of equilibria.}
\end{description}

By \eqref{cis}, the CIS $\mathcal{C}_{d\star}$ is reversible, e.g.~with reversor
\begin{equation}\label{cirev}
\mathcal{R}:\ \quad A_j\, \leftrightarrow\, B_{d-j}
\end{equation}
for $0\leq j\leq d$; see \eqref{rev}.
Analogously to fig.~\ref{ci3}, indeed, the automorphism $\mathcal{R}$ of the CIS $\mathcal{C}_{d\star}$ rotates the pointed connection graph by $180^\circ$, reversing the direction of all edges.
On the boundary sphere $\Sigma^{d-1}$, where all edges actually mean heteroclinic orbits, this establishes time reversibility \eqref{rev} of the actual Chafee-Infante connection graph $\mathcal{C}_d$\,.

In section \ref{Conn}, Morse-shifted versions $\mathrm{CIS}^j_d$ of the CIS  $\mathcal{C}_{d\star}$ will play a central role in our description of the 3-nose pointed connection graphs $\mathcal{C}_{rq\star}$\,; see theorem \ref{connthm}.

\section{Sturm permutations and equilibrium labels}\label{Sigma}

In this section we develop our setting and notation to disentangle the ODE 3-nose meanders, and to distill the PDE connection graphs from them.
Although each of these two viewpoints is aesthetically pleasing, in its own right, their precise relationship is somewhat tricky, notationally, and not at all intuitive.
We start from explicit expressions \eqref{h0h1odd}, \eqref{h0h1even} for the meander permutation $\sigma_{rq}$ and its trivially equivalent inverse $\sigma_{rq}^{-1}$, in proposition \ref{sigmaprop}.
This provides an explicit proof of the trivial equivalence \eqref{rqqr} claimed in theorem \ref{rqqrthm}; see corollary \ref{rqqrcor}.
In definition \ref{h1def}, we introduce labels which enumerate the vertices of the associated meander, along the horizontal axis $h_1$ and, alternatively, along the meander $h_0$.
Proposition \ref{h0prop} then shows consistency of the two alternative labelings.
Based on these labels, the pointed connection graph $\mathcal{C}_{rq\star}$ will be described in the main theorem \ref{connthm} of the next section.

For any dissipative meander $\mathcal{M}$ with $N$ vertices and vertex set $\mathcal{E}$, let $h_\iota: \{1,\ldots,N\} \rightarrow \mathcal{E}$ denote the enumerations of the vertices along the meander and along the horizontal axis, i.e.~for $\iota=0$ and $\iota=1$, respectively.
This generalizes the Sturm case \eqref{hdef} to dissipative meanders which are not necessarily Morse.
Then \eqref{permdef} defines the associated meander permutation $\sigma$ and the inverse $\sigma^{-1}$ as
	\begin{equation}
	\sigma= h_0^{-1} \circ h_1\, \quad \textrm{and} \quad \sigma^{-1}= h_1^{-1} \circ h_0\,.
	\label{perminvdef}
	\end{equation}

\begin{prop}\label{sigmaprop}
On their $N=2(r+1)(q+1)-1$ vertices, the following two relations define the meander permutations  $\sigma=\sigma_{rq}$ and their trivially equivalent inverses $\sigma_{rq}^{-1}$, simultaneously, via \eqref{perminvdef}\,\emph{:}
\begin{align} 
  \label{h0h1odd}
  h_1\Big(2\big((q+1)j+k\big)+1\Big)\ &=\ h_0\Big(2\big(j+(r+1)k\big)+1\Big) \,;\\
  \label{h0h1even}
    h_1\Big(2\big((q+1)(r-j)+(q-k)+1\big)\Big)\ &=\ h_0\Big(2\big(j+(r+1)k\big)\Big) \,.
\end{align}
The relations hold for all $0\leq j \leq r$ and $0\leq k \leq q$, with the exception of the pair $j=k=0$ in \eqref{h0h1even}.
Note how \eqref{h0h1odd} defines the parity preserving meander permutations at odd arguments, from 1 to $N$, whereas \eqref{h0h1even} addresses even arguments, from 2 to $N-1$.
\end{prop}

\begin{proof}
It is sufficient to show that $\sigma^{-1}:=h_1^{-1}\circ h_0$, defined by \eqref{h0h1odd}, \eqref{h0h1even}, is the inverse of the meander permutation $\sigma_{rq}$\,.
We verify this, separately, for upper and lower arcs of the meander $\mathcal{M}=\mathcal{M}_{pq}$\,, with $p=r(q+1)$ as in \eqref{prq}.

The rainbow lower arcs of $\mathcal{M}_{pq}$ are characterized by the invariant sum
\begin{equation}
\label{lowerrainbow}
\sigma_{rq}^{-1}(n)+\sigma_{rq}^{-1}(n+1) = N+2,
\end{equation}
for all even arguments $n=2,\ldots,N-1$.
Here we use that the meander arc from vertex $A=h_0(n)$ to $B=h_0(n+1)$ is a lower arc, if and only if $n$ is even.
Moreover the positions of $A,B$ along the horizontal axis are enumerated by $h_1^{-1}(A), h_1^{-1}(B)$, respectively.
By \eqref{perminvdef}, this establishes the lower rainbow characterization \eqref{lowerrainbow}.

We now verify that $\sigma^{-1}$ defined via \eqref{h0h1odd}, \eqref{h0h1even} satisfies the same lower rainbow characterization \eqref{lowerrainbow}.
We therefore insert \eqref{h0h1even} for the even argument $n=2\big(j+(r+1)k\big)$ on the right, and \eqref{h0h1odd} for odd $n+1=2\big(j+(r+1)k\big)+1$. 
Then $\sigma^{-1}=h_1^{-1}\circ h_0$ from \eqref{perminvdef} implies
\begin{equation}
\label{lowerrainbowpf}
\begin{aligned}
\sigma^{-1}(n)+\sigma^{-1}(n+1) &= 2\big((q+1)(r-j)+(q-k)+1\big) + 2\big((q+1)j+k\big) + 1 =\\
                                                    &=  2\big((q+1)r+q+1\big)+1 = N+2,
\end{aligned}
\end{equation}
as required by \eqref{lowerrainbow}.

Working with upper arcs from $A=h_0(n)$ to $B=h_0(n+1)$, i.e.~for odd $n$, we can analogously verify that the characterizations
\begin{eqnarray}
\label{upperp}
\sigma_{rq}^{-1}(n)+\sigma_{rq}^{-1}(n+1) \ =& 2p+1 &=\  2r(q+1)+1,\\
\label{upperq}
\sigma_{rq}^{-1}(n)+\sigma_{rq}^{-1}(n+1) \ =& 4p+2q+1 &=\ 4r(q+1)+2q.
\end{eqnarray}
also hold for $\sigma^{-1}$.
Here we have to insert odd $n\in \{1,\ldots,2r(q+1)-1\}$ from the right of \eqref{h0h1odd}, and even $n+1$ from the right of \eqref{h0h1even}, to verify the characterization \eqref{upperp} of the left upper $p$-nest, with $p=r(q+1)$.
The analogous characterization \eqref{upperq} of the right upper $q$-nest has to hold for all odd $n\in \{2r(q+1)+1,\ldots,2(r+1)(q+1)-3\}$.
We omit the somewhat tedious details, which are analogous to \eqref{lowerrainbowpf}.

This then verifies that $\sigma$ indeed describes the same meander $\mathcal{M}_{pq}$ as $\sigma_{rq}$ does.
Therefore $\sigma = \sigma_{rq}$\,, and the proposition is proved.
\end{proof}

To tackle claim \eqref{rqqr} concerning the swap $r \leftrightarrow q$, later, we better keep track of $r$ and $q$, notationally.
We first rewrite \eqref{permdef} as
\begin{equation}
\label{hrq}
\begin{aligned}
    \sigma_{rq} \ &=\  (h_0^{rq})^{-1}\circ h_1^{rq},   \\
    \sigma_{qr} \ &=\ (h_0^{qr})^{-1}\circ h_1^{qr} . 
\end{aligned}
\end{equation}
We use the consistent abbreviations
\begin{equation}
\label{nrq}
\begin{aligned}
    n^{rj}_{qk} \ &:=\  2\big((r+1)k+j\big),   \\
    n^{qk}_{rj} \ &:=\  2\big((q+1)j+k\big).
\end{aligned}
\end{equation}
The redundant left subscripts $q,r$ just keep track of the ranges of the right subscripts $k,j$, respectively.
We can then rewrite the defining expressions \eqref{h0h1odd}, \eqref{h0h1even} for $\sigma_{rq}$ and $(\sigma_{rq})^{-1}$ more pedantically, but also more concisely, as
\begin{align} 
  \label{rqodd}
  h_1^{rq}(n^{qk}_{rj}+1)\ &=\ h_0^{rq}(n^{rj}_{qk} +1) \,,\\
  \label{rqeven}
  h_1^{rq}(n^{q,q+1-k}_{r,r-j})\ &=\ h_0^{rq}(n^{rj}_{qk} ) \,.
\end{align}
The analogous defining expressions for  $\sigma_{qr}$ and $(\sigma_{qr})^{-1}$  are simply obtained by the swaps $r\leftrightarrow q$ and $j\leftrightarrow k$ to be
\begin{align} 
  \label{qrodd}
  h_1^{qr}(n^{rj}_{qk}+1)\ &=\ h_0^{qr}(n^{qk}_{rj} +1) \,,\\
  \label{qreven}
  h_1^{qr}(n^{r,r+1-j}_{q,q-k})\ &=\ h_0^{qr}(n^{qk}_{rj} ) \,.
\end{align}

\begin{defi}\label{h1def}
For general $r,q\geq 1$, we label the vertices of the meander permutation $\sigma_{rq}$\,, etc., along the horizontal $h_1$-axis, as\begin{equation}
\label{AB1}
\begin{array}{rcl }
j \ \mathrm{odd:} & & \\[2pt]
     A^{rj}_{qk} & :=  & h_1^{rq}\big((q+1)(j-1)+k+1\big)   \\[2pt]
     B^{rj}_{qk}      & :=  &  h_1^{rq}\big((q+1)j+(q-k)+1\big)   \\[2pt]
     & & \\
j \ \mathrm{even:} & &  \\[2pt]
     A^{rj}_{qk} & :=  & h_1^{rq}\big((q+1)(2r+1-j)+(q-k)+1\big)   \\[2pt]
     B^{rj}_{qk}      & :=  &  h_1^{rq}\big((q+1)(2r-j)+k+1\big)  
\end{array}
\end{equation}
Along the meander path $h_0^{rq}$, alternatively, we label the same vertices as
\begin{equation}
\label{AB0}
\begin{array}{rcl }
k \ \mathrm{odd:} & & \\
     A^{rj}_{qk} & :=  & h_0^{rq}\big((r+1)(2q+1-k)+(2r+1-j)\big)   \\[2pt]
     B^{rj}_{qk}      & :=  &  h_0^{rq}\big((r+1)(2q+1-k)+ j \big)   \\
     & & \\
k \ \mathrm{even:} & &  \\
     A^{rj}_{qk} & :=  & h_0^{rq}\big((r+1)k+j\big)   \\[2pt]
     B^{rj}_{qk}      & :=  &  h_0^{rq}\big((r+1)(k+1)+(r-j)\big)  
\end{array}
\end{equation}
\end{defi}

\begin{figure}[t!]
\begin{center}
\centering \includegraphics[width=\textwidth]{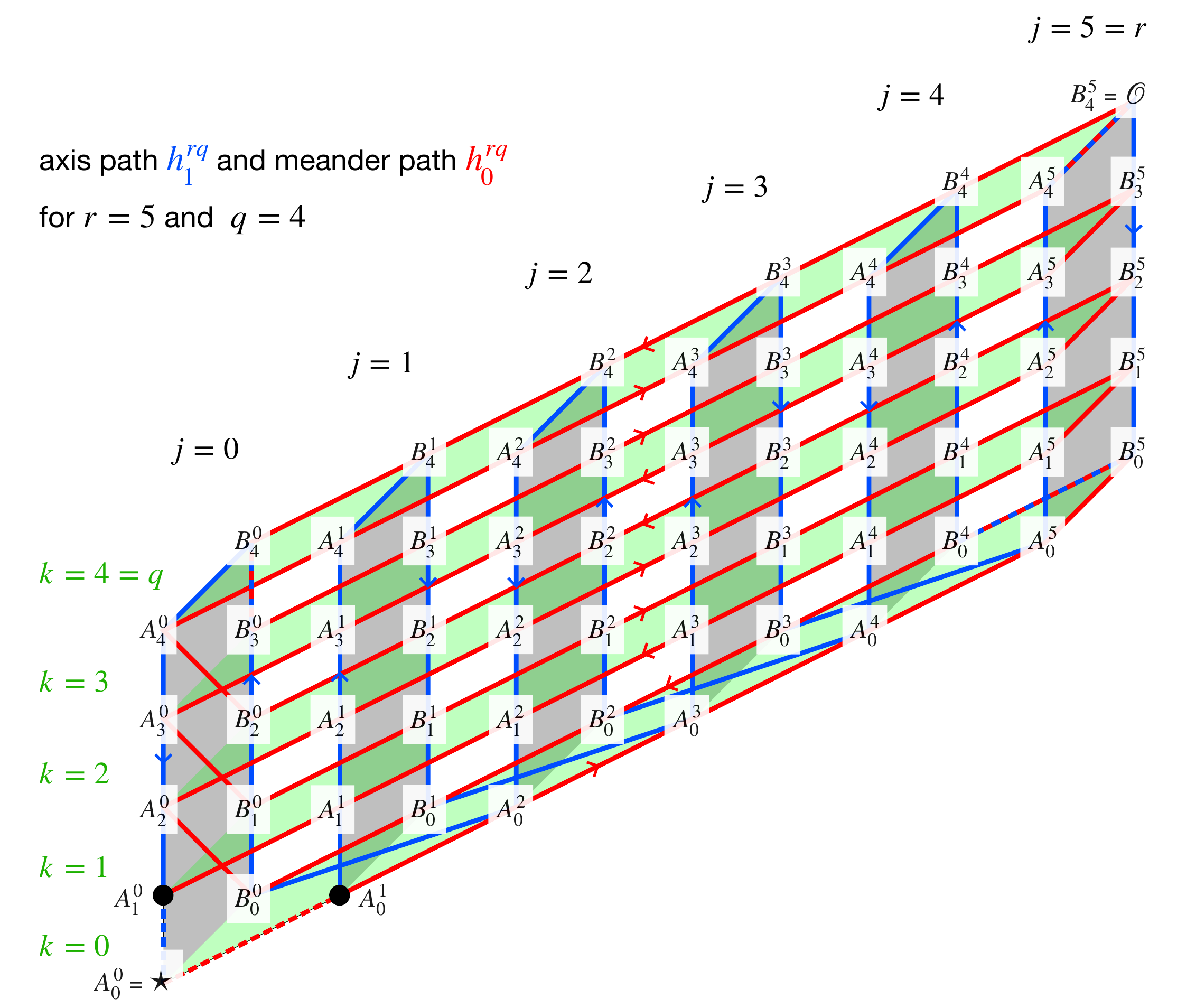}
\caption{\emph{
Visualization of the directed axis path $h^{rq}_1$ (blue) and the meander path $h^{rq}_0$ (red), for $r=5$ and $q=4$.
See definition \ref{h1def}, and in particular \eqref{AB1}, \eqref{AB0}.
We will prove later that horizontal rows indicate equal Morse levels; see \eqref{iLjk} and theorem \ref{connthm}.
Both paths emerge from the $i=0$ sink vertex $A^1_0$ and terminate at sink $A^0_1$ (black dots).
We do not draw the associated meander $\mathcal{M}_{pq}\,,\ p=r(q+1)$, of $N$=$2(r+1)(q+1)-1=59$ vertices.
For a general meander template see fig.~\ref{qnest}(b).
The three meander noses, where the red and blue paths overlap, are dashed blue-red.
In order $h^{rq}_1$ of appearance along the horizontal axis, the three noses are the nose $A^5_4 B^5_4$ of the left upper $p$-nest, the lower rainbow nose $B^5_0 B^4_0$, and the nose $B^0_3 B^0_4$ of the right upper $q$-nest.
Note the distinguished extra vertex $h^{rq}_0(0):= \star = A^0_0 =: h^{rq}_1(N+1)$, prepended to $h^{rq}_0$ by a dashed red line, and appended to $h^{rq}_1$ in dashed blue.
The formal Chafee-Infante sequences $\mathrm{CI}^j$ of \eqref{CI^j} are the boundaries of the gray shaded regions, for $j=0,\ldots,r$.
Shaded green marks the regions $\mathrm{CI}_k\,, \ k=0,\ldots,q$, of \eqref{CI_k}.
}}
\label{h0h1fig}
\end{center} 
\end{figure}

See fig.~\ref{h0h1fig} below for an illustration of the case $r=5,\ q=4$.

\begin{prop}\label{h0prop}
The above vertex labelings \eqref{AB1} and \eqref{AB0} are consistent.
\end{prop}
\begin{proof}
To compare the two definitions, we first substitute $j=2j'+1$ for odd $j$ and $j=2j'$ for even $j$.
Similarly we distinguish the parities of $k=2k',\ 2k'+1$.
Omitting primes, after substitution, and invoking notation \eqref{nrq}, we patiently obtain the following 8-fold path of duplicate definitions:
\begin{equation}
\label{A01}
\begin{array}{lclcl}
     h_1^{rq}(n^{q,q+1-k}_{r,r-j}) & =  & A^{r,2j}_{q,2k} & = & h_0^{rq}(n^{rj}_{qk})   \\[4pt]
     h_1^{rq}(n^{q,q-k}_{r,r-j}+1) & =  & A^{r,2j}_{q,2k+1} & = & h_0^{rq}(n^{r,r-j}_{q,q-k}+1)   \\[4pt]
     h_1^{rq}(n^{qk}_{rj}+1)       & =  & A^{r,2j+1}_{q,2k} & = & h_0^{rq}(n^{rj}_{qk}+1)   \\[4pt]
     h_1^{rq}(n^{q,k+1}_{rj})      & =  & A^{r,2j+1}_{q,2k+1} & = & h_0^{rq}(n^{r,r-j}_{q,q-k})   \\
\end{array}
\end{equation}
\begin{equation}
\label{B01}
\begin{array}{lclcl}
     h_1^{rq}(n^{qk}_{r,r-j}+1) & =  & B^{r,2j}_{q,2k} & = & h_0^{rq}(n^{r,r-j}_{qk}+1)   \\[4pt]
     h_1^{rq}(n^{q,k+1}_{r,r-j}) & =  & B^{r,2j}_{q,2k+1} & = & h_0^{rq}(n^{rj}_{q,q-k})   \\[4pt]
     h_1^{rq}(n^{q,q+1-k}_{rj})       & =  & B^{r,2j+1}_{q,2k} & = & h_0^{rq}(n^{r,r-j}_{qk})   \\[4pt]
     h_1^{rq}(n^{q,q-k}_{rj}+1)      & =  & B^{r,2j+1}_{q,2k+1} & = & h_0^{rq}(n^{rj}_{q,q-k}+1)   \\
\end{array}
\end{equation}
Of course, second subscripts of $L=A,B$ are bound to enumerate the even/odd numbers in $\{0,\ldots,q\}$ here, and second superscripts range in $\{0,\ldots,r\}$.
If we now compare the left and right entries in each row, we recover definition \eqref{h0h1odd}, \eqref{h0h1even} of $\sigma_{rq}$ in the guise of \eqref{rqodd}, \eqref{rqeven}.
Therefore the duplicate definitions \eqref{AB1} and \eqref{AB0} of our vertex labels are consistent.
\end{proof}

After these preliminaries, we can now address the swap $r \leftrightarrow q$.
First consider the meander and axis paths 
\begin{equation}
\label{hrqqrlab}
\begin{aligned}
h_\iota^{rq}:\ \{1,\ldots,N\} \rightarrow \mathcal{E}^r_q\,, \\[2pt]
h_\iota^{qr}:\ \{1,\ldots,N\} \rightarrow \mathcal{E}^q_r\,,
\end{aligned}
\end{equation}
for $\iota=0,1$.
Here the vertex set $\mathcal{E}^r_q$ collects all vertex labels $L^{rj}_{qk}$\,, for tags $L \in \{A,B\}$ and all $j,k$.
Similarly, $\mathcal{E}^q_r$ collects all $L^{qk}_{rj}$\,.
We relate the two vertex sets by the label map
\begin{equation}
\label{L}
\begin{aligned}
\Lambda:\quad \mathcal{E}^q_r &\rightarrow \mathcal{E}^r_q\,, \\[2pt]
L^{qk}_{rj} &\mapsto L^{rj}_{qk}\,.
\end{aligned}
\end{equation}
Only for $r=q$, the label map $\Lambda$ is a selfmap of $\mathcal{E}^q_r = \mathcal{E}^r_q$ and, in fact, and involution.
The label map $\Lambda$ intertwines the four axis and meander paths $h^{qr}_\iota$ and $h^{rq}_\iota$ as follows.

\begin{lem}\label{rqqrlem}
With the above notation and, in particular, the label map \eqref{L} and the rotation $\kappa$ from \eqref{flip}, the following holds true for all $r,q\geq 1$ and $\iota=0,1$\,\emph{:}
\begin{equation}
\label{hrqqr}
\Lambda \,\circ\, h^{qr}_\iota = h^{rq}_{1-\iota}\,\circ\, \kappa\,.
\end{equation}
\end{lem}

\begin{proof}
We check claim \eqref{hrqqr} via the explicit 8-fold path \eqref{A01}, \eqref{B01}, for the particular arguments $n=n^{rj}_{qk}+1$ and $\iota=1$.
The remaining cases are similar.

The swap-invariant action \eqref{flip} of $\kappa$ with $N+1=2(r+1)(q+1)$ implies
\begin{align}
  \label{koddrq}
  \kappa(n^{rj}_{qk}+1)&=n^{r,r-j}_{q,q-k}+1   \,,\\
    \label{koddqr}
  \kappa(n^{qk}_{rj}+1)&=n^{q,q-k}_{r,r-j}+1   \,.
\end{align}
Here we only have to check \eqref{koddrq}, explicitly, because \eqref{koddqr} follows by swap.

We can now use \eqref{koddrq}, and the right hand side of the second line of \eqref{A01}, to evaluate the right hand side of \eqref{hrqqr}:
\begin{equation}
\label{hrqpf}
(h^{rq}_{1-\iota}\,\circ\kappa)\; (n^{rj}_{qk}+1) \ =\   h^{rq}_0(n^{r,r-j}_{q,q-k}+1)=A^{r,2j}_{q,2k+1} \,.
\end{equation}
To evaluate the left hand side of \eqref{hrqqr}, we invoke the left hand side of the third line of \eqref{A01} -- albeit with switched roles of $r$ and $q$, as well as $j$ and $k$:
\begin{equation}
\label{hqrpf}
(\Lambda \,\circ\, h^{qr}_\iota)\; (n^{rj}_{qk}+1) \ =\  \Lambda\big( h^{qr}_1(n^{rj}_{qk}+1) \big) = \Lambda(A^{q,2k+1}_{r,2j})\,.
\end{equation}
In view of the label map \eqref{L}, the two expressions \eqref{hrqpf} and \eqref{hqrpf} coincide.
This proves the lemma.
\end{proof}

\begin{cor}\label{rqqrcor}
Claim \eqref{rqqr} of theorem \ref{rqqrthm} holds true, for all $r,q\geq1$.
In the Morse case, the label map $\Lambda$ of \eqref{L} also expresses a trivial equivalence of connection graphs:
\begin{equation}
\label{v1v2rqqr}
  \big(E_1 \leadsto_{rq} E_2 \big) \quad \Leftrightarrow \quad \big(\Lambda^{-1} E_1 \leadsto_{qr} \Lambda^{-1} E_2\big)\,.
\end{equation}
Here $E_1,E_2 \in \mathcal{E}^r_q$ refer to equilibrium labels in $ \mathcal{E}^r_q$\,, whereas $\Lambda^{-1}E_1,\Lambda^{-1}E_2 \in \mathcal{E}^q_r$\,.
The connection symbol $\leadsto_{rq}$ on the left refers to the Sturm permutation $\sigma_{rq}$\,, whereas $\leadsto_{qr}$ on the right refers to $\sigma_{qr}$\,. In other words, the label map $\Lambda$ extends to a direction preserving isomorphism of the connection di-graphs with their respective labelings:
\begin{equation}
\label{Lqrrq}
\Lambda \mathcal{C}_{qr} \cong \mathcal{C}_{rq}\,.
\end{equation}
\end{cor}

\begin{proof}
We insert \eqref{hrqqr} in the definitions \eqref{hrq} of $\sigma_{qr}$ and $\sigma_{rq}$ to obtain 
\begin{equation}
\label{rqqrpf}
\sigma_{qr} = (h_0^{qr})^{-1}\circ h_1^{qr} = (\Lambda h_0^{qr})^{-1}\circ (\Lambda h_1^{qr}) = (h_1^{rq} \kappa )^{-1}\circ (h_0^{rq} \kappa) = \kappa \sigma_{rq}^{-1} \kappa\,.
\end{equation}
This proves claim \eqref{rqqr}.

To prove claim \eqref{v1v2rqqr}, we first introduce the more detailed notation
\begin{equation}
\label{(h0h1)v1v2}
(h_0,h_1):\  E_1 \leadsto E_2
\end{equation}
to indicate a heteroclinic orbit under $\sigma := h_0^{-1}\circ h_1$.
For $\nu=1,2$, each equilibrium  $v_\nu(x)$ is labeled by the \emph{same} label $E_\nu$ under both labeling paths $h_\iota: \{1,\ldots,N\} \rightarrow \mathcal{E},\ \iota=0,1$.
This meticulous notation is required to keep track of the \emph{different} vertex sets $\mathcal{E}^r_q$ and  $\mathcal{E}^q_r$ which we use, eventually, to describe the connection graphs  $\mathcal{C}_{rq}$ and $\mathcal{C}_{qr}$ of $\sigma_{rq}$ and $\sigma_{qr}$\,, respectively.

Proceeding with care, we begin with the trivial equivalence
\begin{equation}
\label{rotconn}
   \big( (h_0,h_1):\ E_1\, \leadsto_\sigma\, E_2 \big) \quad \Leftrightarrow \quad \big( (h_0\kappa,h_1\kappa):\ E_1 \,\leadsto_{\kappa\sigma\kappa}\, E_2 \big)\,,
\end{equation}
for any meander permutation $\sigma=h_0^{-1}\circ h_1$ and its trivial equivalent $\kappa\sigma\kappa$.
Note how \emph{the same vertex labels} $E_1$ and $E_2$ appear on the right and on the left of the equivalence.
Indeed, the label $E_\nu$ of an equilibrium profile $v_\nu(x)$, on the left, has to be replaced by the label of the equilibrium profile $\kappa v_\nu(x)=-v_\nu(x)$, on the right; see \eqref{rot}.
The proper \emph{vertex label} of that equilibrium $-v_\nu(x)$, however, remains the same $E_\nu$, on the right and on the left.
The reason is that $\kappa$, just like the minus sign, reverses the order of enumeration on either boundary; see \eqref{flip}, \eqref{rotsig}.

For the inverse permutation $\sigma^{-1}=h_1^{-1}\circ h_0$ we similarly obtain
\begin{equation}
\label{invconn}
   \big( (h_0,h_1):\ E_1 \,\leadsto_\sigma\, E_2 \big) \quad \Leftrightarrow \quad \big( (h_1,h_0):\ E_1 \,\leadsto_{\sigma^{-1}}\, E_2 \big)\,.
\end{equation}
Indeed, this time the same vertex label $E_\nu$ of an equilibrium profile $v_\nu(x)$, on the left, applies to the equilibrium profile $v_\nu(1-x)$, on the right;
see \eqref{flip}, \eqref{invsig}.

For $h_\iota:=h^{rq}_\iota$\,, we can now combine \eqref{rotconn} and \eqref{invconn}, successively, with \eqref{rqqr} and \eqref{hrqqr} to obtain
\begin{equation}
\label{connrqqrpf}
\begin{array}{ccrccc}
  E_1 \leadsto_{rq} E_2  & \Leftrightarrow& (h^{rq}_0,h^{rq}_1):& E_1 \leadsto E_2  \quad &\Leftrightarrow & \\[3pt]
&\Leftrightarrow  & (h^{rq}_0\kappa,h^{rq}_1\kappa):& E_1 \leadsto E_2  \quad &\Leftrightarrow &\\[3pt]
&\Leftrightarrow  & (h^{rq}_1\kappa,h^{rq}_0\kappa):& E_1 \leadsto E_2  \quad &\Leftrightarrow &\\[3pt]
&\Leftrightarrow  &  (\Lambda h^{qr}_0, \Lambda h^{qr}_1):& E_1 \leadsto E_2  \quad &\Leftrightarrow &\\[3pt]
&\Leftrightarrow  & (h^{qr}_0, h^{qr}_1):& \Lambda^{-1} E_1 \leadsto \Lambda^{-1} E_2  \quad &\Leftrightarrow & \Lambda^{-1} E_1 \leadsto_{qr} \Lambda^{-1} E_2 
\end{array}
\end{equation}
Note how $\Lambda^{-1}$ has relabeled vertices $E_\nu$ in the last line, only.
This proves the corollary.
\end{proof}

See fig.~\ref{h0h1fig} again, for an illustration of the case $r=5,\ q=4$.

It is instructive to compare the succession of labels along the paths $h^{rq}_1$ and $h^{rq}_0$.
For $h^{rq}_1$, tags $L\in\{A,B\}$, and omitting $r,q$ again, we use the abbreviations
\begin{equation}
\label{h1seq}
L^j_\nearrow := L^j_0 \ldots L^j_q\,,  \qquad  L^j_\searrow := L^j_q \ldots L^j_0\,. 
\end{equation}
For further abbreviation, and as a prequel to the Chafee-Infante stacks $\mathrm{CIS}^j$, we also define formal Chafee-Infante sequences $\mathrm{CI}^j$ as
\begin{equation}
\label{CI^j}
\mathrm{CI}^j := 
\begin{cases}
      A^j_\nearrow\ B^j_\searrow\,, & \text{for } j \text{ odd}, \\[2pt]
      B^j_\nearrow\ A^j_\searrow\,, & \text{for } j \text{ even}.
\end{cases}
\end{equation}
For odd $r$, \eqref{AB1} then determines the vertex labels along $h^{rq}_1$ to be
\begin{equation}
\label{h1labels-oddr}
\begin{array}{lcccccccc}
      h^{rq}_1: & A^1_\nearrow\ B^1_\searrow\ & A^3_\nearrow\ B^3_\searrow\  & \ldots \ 
      &A^r_\nearrow\ B^r_\searrow \ & B^{r-1}_\nearrow\ A^{r-1}_\searrow\ & \ldots \ 
      & B^2_\nearrow\ A^2_\searrow\ &B^0_\nearrow\ A^0_\searrow 
      \\[2pt]
      & \mathrm{CI}^1 & \mathrm{CI}^3 & \ldots
      &\mathrm{CI}^r  &\mathrm{CI}^{r-1} & \ldots
       & \mathrm{CI}^2    &\mathrm{CI}^0  
\end{array}
\end{equation}
Here the superscripts are odd, in the left part, ascending from 1 to $r$.
In the right part, the even superscripts descend from $r-1$ to 0.

For even $r$, we obtain analogously
\begin{equation}
\label{h1labels-evenr}
\begin{array}{lcccccccc}
      h^{rq}_1: & A^1_\nearrow\ B^1_\searrow\ & A^3_\nearrow\ B^3_\searrow\  & \ldots \ 
      &A^{r-1}_\nearrow\ B^{r-1}_\searrow \ & B^r_\nearrow\ A^r_\searrow\ & \ldots \ 
      & B^2_\nearrow\ A^2_\searrow\ &B^0_\nearrow\ A^0_\searrow 
      \\[2pt]
      & \mathrm{CI}^1 & \mathrm{CI}^3 & \ldots
      &\mathrm{CI}^{r-1}  &\mathrm{CI}^r & \ldots
       & \mathrm{CI}^2    &\mathrm{CI}^0  
\end{array}
\end{equation}
Note how we have tacitly appended the distinguished vertex $h^{rq}_1(N+1):= \star = A^0_0$ at the very end $N+1=2(r+1)(q+1)$ of the enumeration, in either case.

For $h^{rq}_0$ we use analogous abbreviations
\begin{align}
\label{h0seq}
L^\nearrow_k &:= L^0_k \ldots L^r_k\,,  \qquad  L^\searrow_k := L^r_k \ldots L^0_k\,, \\[4pt]
\label{CI_k}
\mathrm{CI}_k &:= 
\begin{cases}
      A^\nearrow_k \ B^\searrow_k\,, & \text{for } k \text{ odd}, \\[\medskipamount]
      B^\nearrow_k\ A^\searrow_k \,, & \text{for } k \text{ even}.
\end{cases}
\end{align}
For odd $q$, \eqref{AB0} then determines the vertex labels along $h^{rq}_0$ to be
\begin{equation}
\label{h0labels-oddq}
\begin{array}{lcccccccc}
      h^{rq}_0: & A^\nearrow_0\ B^\searrow_0\  & A^\nearrow_2 \  B^\searrow_2 \  & \ldots \ 
      &A^\nearrow_{q-1} \ B^\searrow_{q-1} \ & B^\nearrow_q \ A^\searrow_q \ & \ldots \ 
      & B^\nearrow_3 \ A^\searrow_3 \ &B^\nearrow_1 \ A^\searrow _1
      \\[2pt]
      & \mathrm{CI}_0 & \mathrm{CI}_2 & \ldots
      &\mathrm{CI}_{q-1}  &\mathrm{CI}_q & \ldots
       & \mathrm{CI}_3    &\mathrm{CI}_1  
\end{array}
\end{equation}
This time, the subscripts are even, in the left part, and ascending from 0 to $q-1$.
In the right part, the odd subscripts descend from $q$ to 1.

For even $q$, we obtain analogously
\begin{equation}
\label{h0labels-evenq}
\begin{array}{lcccccccc}
      h^{rq}_0: & A^\nearrow_0\ B^\searrow_0\  & A^\nearrow_2 \  B^\searrow_2 \  & \ldots \ 
      &A^\nearrow_q \ B^\searrow_q \ & B^\nearrow_{q-1} \ A^\searrow_{q-1} \ & \ldots \ 
      & B^\nearrow_3 \ A^\searrow_3 \ &B^\nearrow_1 \ A^\searrow _1
      \\[2pt]
      & \mathrm{CI}_0 & \mathrm{CI}_2 & \ldots
      &\mathrm{CI}_q  &\mathrm{CI}_{q-1} & \ldots
       & \mathrm{CI}_3    &\mathrm{CI}_1  
\end{array}
\end{equation}
This time, the distinguished vertex $h^{rq}_0(0):= \star = A^0_0$ has been prepended at the very beginning of the enumeration, in either case.

In summary, \eqref{h1seq}--\eqref{h0labels-evenq} enumerate the vertex paths
\begin{equation}
\label{h0h1labels}
\begin{array}{lccccccc}
      h^{rq}_0: & \mathrm{CI}_0 & \mathrm{CI}_2 & \ldots
      &\mathrm{CI}_k  & \ldots
       & \mathrm{CI}_3    &\mathrm{CI}_1  
      \\[2pt]
      h^{rq}_1: & \mathrm{CI}^1 & \mathrm{CI}^3 & \ldots
        &\mathrm{CI}^j & \ldots
       & \mathrm{CI}^2    &\mathrm{CI}^0  
\end{array}
\end{equation}
for $0\leq j\leq r$ and $0\leq k\leq q$ with the indicated succession of parities.

\section{Chafee-Infante lattices}\label{Conn}

In this section we address the remaining four theorems \ref{krmuthm}, \ref{rqqrthm}, \ref{krballthm}, and \ref{krrevthm}, of section \ref{Intro}, on the primitive 3-nose Sturm attractors $\mathcal{A}_{rq}$\,, their dissipative Morse meanders $\mathcal{M}_{r(q+1),q}$\,, and their entourage of Sturm permutations $\sigma_{rq}$ and connection graphs $\mathcal{C}_{rq}$\,.
Theorem \ref{connthm} below describes the pointed connection graphs $\mathcal{C}_{rq\star}$\,, based on two alternative decompositions $\mathrm{CIS}^{rq}$ and $\mathrm{CIS}_{rq}$ into Morse-shifted Chafee-Infante stacks (CIS).
In sections \ref{krmupf}--\ref{krrevpf} we then show how the connection graph of theorem \ref{connthm} implies the four remaining theorems.
The proof of theorem \ref{connthm} itself will be postponed to section \ref{Rec}.

Consider the vertices $L^j_k$\,, with tags $L\in \{A,B\}$, in the Chafee-Infante sequence $\mathrm{CI}^j$ introduced in \eqref{CI^j}.
The labels fit the definition \eqref{cis} of a Chafee-Infante stack:
for any fixed $j \in \{0, \ldots, r\}$, the obvious vertex bijection $\mathcal{C}_{q\star} \rightarrow \mathrm{CI}^j$ is
\begin{equation}
\label{toCI^j}
L_k \mapsto L^j_k\,,
\end{equation}
for all $0\leq k\leq q$. 
We define the \emph{vertical Morse-shifted Chafee-Infante stack $\mathrm{CIS}^j$ of height $q+1$}, on the vertices provided by $\mathrm{CI}^j$, by inducing all directed edges on $\mathrm{CIS}^j$ from the Chafee-Infante stack $\mathcal{C}_{q\star}$\,, via the map \ref{toCI^j}.
To equip the graded graph $\mathrm{CIS}^j$ with a ranking function, analogously to \eqref{iLj}, we define formal Morse levels of any vertex $L^j_k \in \mathrm{CIS}^j$ as
\begin{equation}
\label{iLjk}
i(L^j_k)\ =\ \begin{cases}
      \ j+k-1& \text{for } L=A, \\
      \ j+k& \text{for } L=B.
\end{cases}
\end{equation}
In other words, the Morse levels of $\mathrm{CIS}^j$ have been shifted by the superscript $j$, compared to the standard pointed Chafee-Infante connection graph $\mathcal{C}_q$\,.
Note how all stacks $\mathrm{CIS}^j$ are mutually vertex disjoint, by construction.

For any fixed $k \in \{0, \ldots, q\}$, we similarly lift $\mathcal{C}_{r\star} \rightarrow \mathrm{CI}_k$ by
\begin{equation}
\label{toCI_k}
L_j \mapsto L^j_k\,,
\end{equation}
for all $0\leq j \leq r$. 
Analogously to $\mathrm{CIS}^j$ above, this defines the mutually disjoint \emph{slanted Morse-shifted Chafee-Infante stacks} $\mathrm{CIS}_k$ of height $r+1$.

\begin{figure}[t!]
\begin{center}
\centering \includegraphics[width=\textwidth]{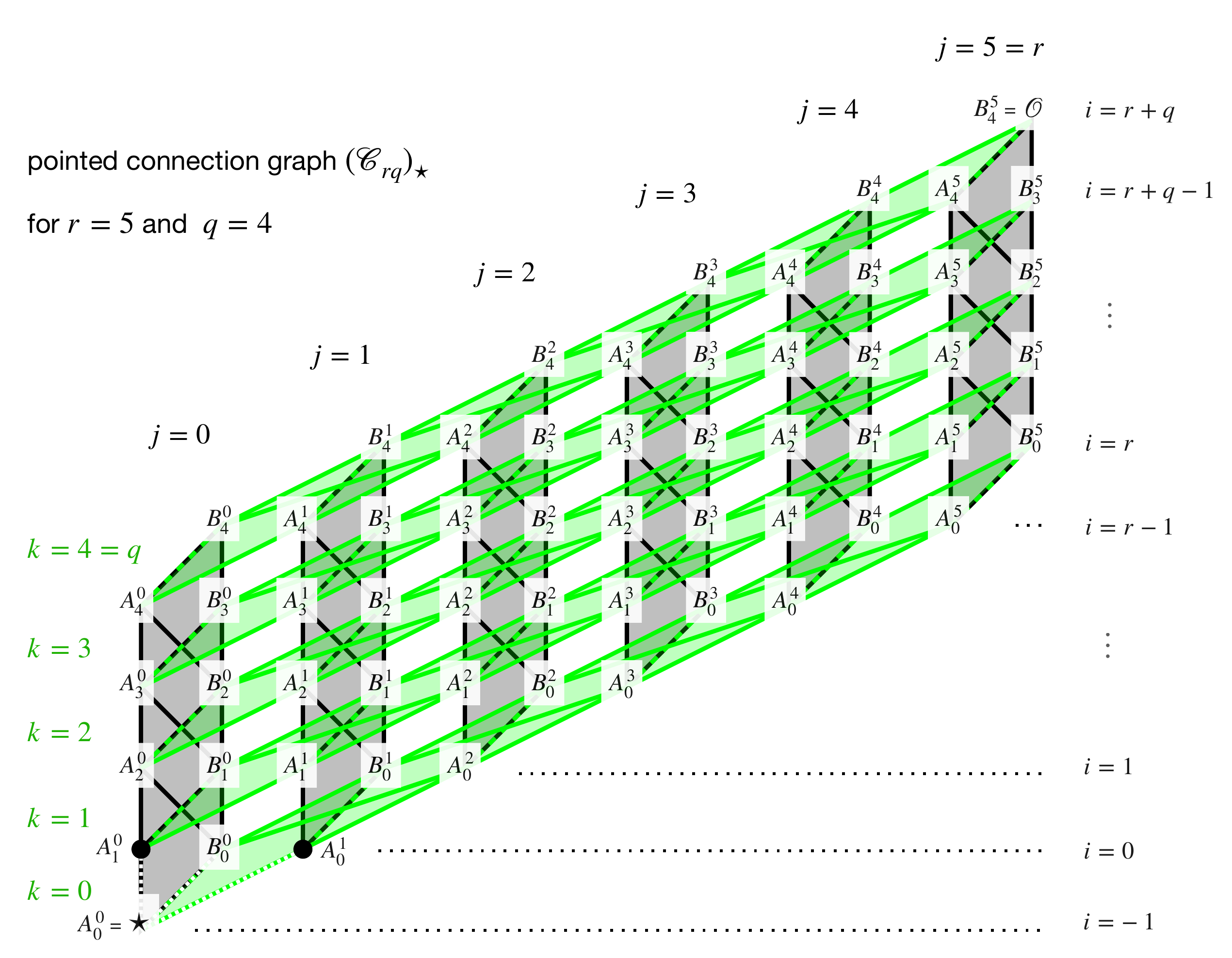}
\caption{\emph{
The Chafee-Infante lattice of the pointed connection graph $\mathcal{C}_{rq\star}$\,, for $r=5$ and $q=4$.
See \eqref{toCI^j}--\eqref{CIS_rq}, and theorem \ref{connthm}.
Compare \eqref{AB1}, \eqref{AB0}, and fig.~\ref{h0h1fig} for labels, annotations, and gray/green shading of Chafee-Infante sequences, which become the two stack decompositions \eqref{CIS^rq} and \eqref{CIS_rq} of $\mathcal{C}_{rq\star}$ into Morse-shifted Chafee-Infante stacks $\mathrm{CIS}^j$ and $\mathrm{CIS}_k$\,.
The graph is graded by the ranking function $i$ of Morse levels, as indicated on the right.
Morse-adjacent heteroclinic orbits in gray shaded vertical Chafee-Infante stacks $\mathrm{CIS}^j$ are black edges.
Green edges belong to green shaded slanted stacks $\mathrm{CIS}_k$\,.
Shared edges of the two stack decompositions are dashed green-black.
Downwards arrows have been omitted.
Lower left: the distinguished vertex $\star=A^0_0$ resides at formal Morse level $i=-1$.
The dotted edges towards the distinguished vertex $\star$ in the stacks $\mathrm{CIS}^0$ and $\mathrm{CIS}_0$ are formal, and do not correspond to actual heteroclinic orbits.\\
Note how the connection graph strictly contains both blue/red paths $h^{rq}_\iota$ of fig.~\ref{h0h1fig}.
Indeed, path-adjacent vertices are Morse-adjacent and connected; see \eqref{i} and \eqref{wolfrum}, \eqref{block}.
}}
\label{connfig}
\end{center} 
\end{figure}

We combine the above Chafee-Infante stacks, respectively, into two graded di-graphs
\begin{align}
\label{CIS^rq}
\mathrm{CIS}^{rq} \ &:=\ \bigcup_{j=0}^r \mathrm{CIS}^j \quad \textrm{(vertical)},\\
\label{CIS_rq}
\mathrm{CIS}_{rq} \ &:=\ \bigcup_{k=0}^q \mathrm{CIS}_k  \quad \textrm{(slanted)}.
\end{align}
By construction, each stack decomposition defines a graded graph on the same $2(r+1)(q+1)$ vertices as the pointed connection graph  $\mathcal{C}_{rq\star}$ which we aim to describe.
The ranking function by Morse levels $i$ will turn out to be the same, for both graded graphs.
In fact, $\mathrm{CIS}^{rq}$ decomposes the vertices of $\mathcal{C}_{rq\star}$ into the $r+1$ Morse-shifted vertical Chafee-Infante stacks $\mathrm{CIS}^j$, each of height $q+1$ with vertices $L^{rj}_{qk},\ k=0,\ldots,q$, and tags $L=A,B$.
Alternatively, $\mathrm{CIS}_{rq}$ rearranges the same vertices $L^{rj}_{qk}$ into the $q+1$ Morse-shifted slanted Chafee-Infante stacks $\mathrm{CIS}_k$ of height $r+1$.

\begin{thm}\label{connthm}
Let $rq >1$. 
In the above notation, the pointed connection graph $\mathcal{C}_{rq\star}$ is the union of the two stack decompositions \eqref{CIS^rq} and \eqref{CIS_rq}, i.e.
\begin{equation}
\label{lattice}
\mathcal{C}_{rq\star} \ = \ \mathrm{CIS}^{rq} \,\cup \mathrm{CIS}_{rq}\,.
\end{equation}
In particular, the sums $j+k$ in the formal Morse levels \eqref{iLjk} indicate actual Morse indices.
\end{thm}

See fig.~\ref{connfig} for an illustration of the case $r=5,\ q=4$ with $2(r+1)(q+1)=60$ vertices.
Note the two stack decompositions into $r+1=6$ vertical Chafee-Infante stacks $\mathrm{CIS}^j$ of height $q+1=5$ (shaded gray), and into $q+1=5$ slanted Chafee-Infante stacks $\mathrm{CIS}_k$ of height $r+1=6$ (shaded green).
For an explicit version of the simplest case $r=1,\ q=2$, rotated by the trivial equivalence $\kappa$, see also $\sigma=\sigma_{12}^\kappa,\ \mathcal{A}=\mathcal{A}_{12}^\kappa,\ \mathcal{C}_\star=\mathcal{C}^\kappa_{12\star}$ sketched in fig.~\ref{3ball}(a),(c),(d).
In that case of $2(r+1)(q+1)=12$ vertices, we observe only $r+1=2$ vertical Chafee-Infante stacks $\mathrm{CIS}^j$ of height $q+1=3$ (black), but  3 slanted Chafee-Infante stacks $\mathrm{CIS}_k$ of height 2 (green).
The trivial equivalence $\kappa$ of \eqref{rot}, \eqref{flip} leaves the pointed connection graph $\mathcal{C}^\kappa_{12\star}$ invariant, down to identical labels.

By a slight abuse of terminology, we call the union \eqref{lattice} a \emph{Chafee-Infante lattice}.
Each Chafee-Infante stack $\mathrm{CIS}^j,\ \mathrm{CIS}_k$, like any acyclic digraph, defines a partial order on the vertices.
Even for a single stack, however, this partial order is not a lattice. 
Indeed, the set $\{A^j_{k+1}\,, B^j_k\}$ at Morse level $i=j+k$ possesses two candidates for a ``supremum'', and two candidates for an ``infimum'', given by subscripts raised and lowered by 1, respectively.
In the standard terminology of lattices, this violates their required uniqueness.
On the other hand, the resulting graph in fig.~\ref{connfig} just looks like a lattice fence.
In fact, we could enforce strict compliance with the standard definition, e.g., by an identification of the two tags $L=A,B$ for each $L^{rj}_{kq}$\,.
After such identification, we even obtain a skew $(r+1)\times (q+1)$ grid.
Averse to further petty pedantry, we call \eqref{lattice} a lattice, anyway.

We defer the proof of our main theorem \ref{connthm} to section \ref{Rec}.
Based on theorem \ref{connthm}, for now, we can already prove the remaining four theorems \ref{krmuthm}, \ref{rqqrthm}, \ref{krballthm}, and \ref{krrevthm}.
See also fig.~\ref{connfig}.

\subsection{Proof of theorems \ref{krmuthm} and \ref{rqqrthm}.}\label{krmupf}
The pointed connection graph $\mathcal{C}_{rq\star}$ of theorem \ref{connthm} allows us to count  equilibria via the stack decomposition \eqref{CIS^rq} into the $r+1$ vertical stacks $\mathrm{CIS}^j$ or, alternatively but equivalently, via the stack decomposition \eqref{CIS_rq} into the $q+1$ slanted stacks $\mathrm{CIS}_k$\,.
Therefore we may assume $r\leq q$, without loss of generality, i.e. $\min \{r,q\} = r$ and $\max\{r,q\} = q$.
Since the vertical Chafee-Infante stacks $\mathrm{CIS}^j$ in $\mathcal{C}_{rq\star}$\,, each of height $q$, are shifted upwards in Morse levels by shifts $j=0,\ldots,r$,  their nonzero Morse counts add up to
\begin{equation}\label{krmu'}
\mu_{rq}(i)\ =\ 
\begin{cases}
      3+2i,& \textrm{for}\quad -1\leq i< r;\\
      2(r+1),& \textrm{for}\qquad\  r\leq i < q;\\
      2(r+q)+1-2i,& \textrm{for}\qquad\  q\leq i\leq r+q.
\end{cases}
\end{equation}
For $r\leq q$, this proves the Morse counts $\mu_{rq}(i)$ of claim \eqref{krmu}, and theorem \ref{krmuthm}. 

In particular, the dissipative meander permutations $\sigma_{rq}$ are Morse, hence Sturm.
Since claim \eqref{rqqr} has already been established in corollary \ref{rqqrcor}, this also completes the proof of theorem \ref{rqqrthm}.
\hfill $\bowtie$

\subsection{Proof of theorem \ref{krballthm}.}\label{krballpf}
By the Schoenflies theorem \cite{firo15}, it is sufficient to prove that the single equilibrium $\mathcal{O}=B^r_q$ of top Morse index $i(\mathcal{O})=r+q=\dim \mathcal{A}_{rq}$ connects heteroclinically to all other equilibria $E$. 
In symbols, $B^r_q\leadsto E$.

By Morse-Smale transitivity of the directed edge relation $\leadsto$, it is sufficient to establish a di-path from $B^r_q$ to any equilibrium $E$, in the connection di-graph $\mathcal{C}_{rq}$\,.
This is obvious from the two Chafee-Infante stack decompositions of the pointed connection di-graph $\mathcal{C}_{rq\star} \supset \mathcal{C}_{rq}$ in theorem \ref{connthm}.
Indeed, consider the stack decomposition \eqref{CIS^rq} into vertical stacks $\mathrm{CIS}^j, \ j=0\ldots r$, of height $q+1$ and with top vertices $B^j_q$\,, respectively.
Then the target $E$ is contained in  $\mathrm{CIS}^j$, for some $j$. 
Any top vertex $B^j_q$, on the other hand, is contained in the heteroclinic di-path
\begin{equation}\label{top descent}
\mathcal{O}=B^r_q \leadsto \ldots \leadsto  B^0_q
\end{equation}
asserted by the top level slanted stack $\mathrm{CIS}_q$ of height $r+1$, in stack decomposition \eqref{CIS_rq}.
This proves theorem \ref{krballthm}.
\hfill $\bowtie$

\subsection{Proof of theorem \ref{krrevthm}.}\label{krrevpf}
For the pointed connection graph $\mathcal{C}_{rq\star}$ of theorem \ref{connthm}, fig.~\ref{connfig}, consider the vertex map
\begin{equation}\label{krrev}
\mathcal{R}:\ \quad A^j_k\, \leftrightarrow\, B^{r-j}_{q-k}\,.
\end{equation}
In other words, the involution $\mathcal{R}$ rotates $\mathcal{C}_{rq\star}$ by $180^\circ$.
This maps each of the two stack decompositions \eqref{CIS^rq}, \eqref{CIS_rq} onto itself.
More precisely, the superscripts $j=0,\ldots,r$ of the decomposition into vertical stacks $\mathrm{CIS}^j$ are swapped with the superscripts $r-j$ of the same decomposition.
Similarly, the subscripts $k=0,\ldots,q$ of the slanted stacks $\mathrm{CIS}_k$ are swapped with $q-k$.
Also, the rotation $\mathcal{R}$ swaps the vertex tags $L \in \{A,B\}$.
Therefore, \eqref{krrev} defines an automorphism $\mathcal{R}$ of the pointed connection graph $\mathcal{C}_{rq\star}$\,, mapping vertices to vertices and edges to edges.

It remains to show that $\mathcal{R}$ reverses edge orientation.
We already know that $\mathcal{R}$ swaps subscripts $k\leftrightarrow q-k$ and superscripts $j\leftrightarrow r-j$ of vertices $L^j_k$\,, as well as tags $L \in \{A,B\}$.
By \eqref{iLjk}, this also swaps Morse levels $i=j+k-1$ of $A^j_k$ with $(r-j)+(q-k)=r+q-1-i$ of $B^{r-j}_{q-k}$\,, for $i=-1,\ldots,r+q$.
Since all edges $\leadsto$ in $\mathcal{C}_{rq\star}$ are directed downwards, between adjacent Morse levels $i$, this shows that the reversor $\mathcal{R}$ indeed reverses all arrows -- which proves theorem \ref{krrevthm}. 
 \hfill $\bowtie$

\section{Recursion of meanders}\label{Rec}

This section consists of the proof of our main theorem \ref{connthm} on the structure of the connection graph $\mathcal{C}_{rq}$\,.
We proceed by induction on, both, $r\geq 2$ and $q\geq 2$.
Specifically, our \emph{induction hypothesis} is: 
\begin{equation}
\label{ind}
\begin{aligned}
\textrm{Theorem \ref{connthm} holds for,} &\textrm{ both, the connection graphs} \\ 
\mathcal{C}_{r-1,q}\  &\textrm{ and } \ \mathcal{C}_{r,q-1}\,.
\end{aligned}
\end{equation}
Let $\mathcal{C}$ denote what we \emph{claim} the connection graph $\mathcal{C}_{rq}$ to be, i.e.
\begin{equation}
\label{connclaim}
\mathcal{C} := \mathrm{CIS}^{rq} \cup \mathrm{CIS}_{rq} \setminus \{\star\};
\end{equation}
see fig.~\ref{connfig}

\begin{figure}[t]
\begin{center}
\centering \includegraphics[width=\textwidth]{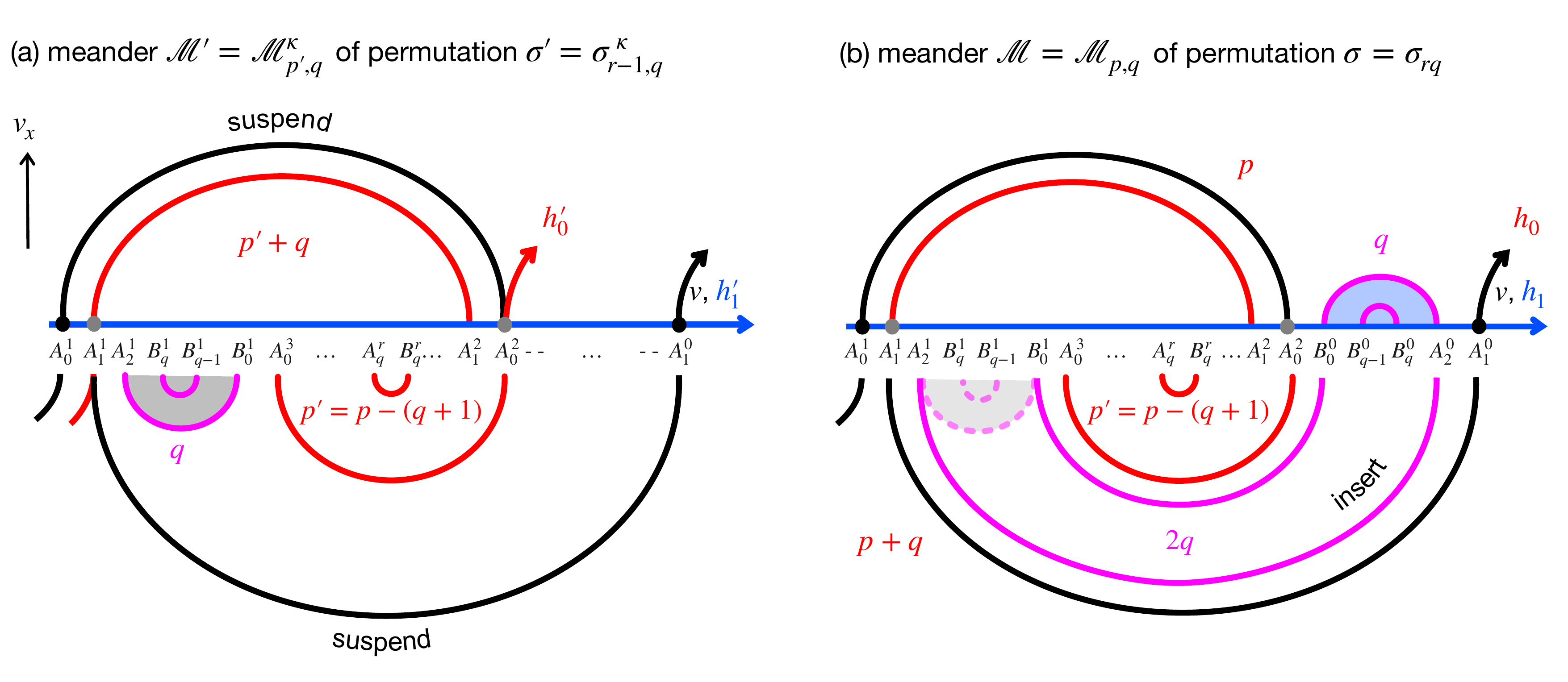}
\caption{\emph{
Induction step from $(r-1,q)$ to $(r,q)$, on the level of meander templates $\mathcal{M}' := \mathcal{M}^{\kappa}_{p'q}$ and $\mathcal{M} := \mathcal{M}_{pq}$\,. 
Here $p' := (r-1)(q+1)\,,\ p := r(q+1)$. \\
(a) Suspension (black) of $\mathcal{M}'$ (red and purple) raises all Morse numbers by 1; see proposition \ref{suspprop}(iv).
Vertex annotation is chosen to coincide with meander $\mathcal{M}$ in (b), for those vertices already present in meander $\mathcal{M}'$.
The $i=0$ start and termination vertices $A^1_1$  and $A^2_0$ of $\mathcal{M}'$ are marked gray.
Their Morse indices become $i=1$, after suspension.
The $i=0$ start and termination vertices $A^1_0$ and $A^0_1$ of the suspension $\widetilde{\mathcal{M}'}$ of $\mathcal{M}'$ are marked black. 
Suspension reverses the direction of path $h_0'$ but not of $h_1'$.\\
(b) Insertion of the $q$-nest in the suspension $\widetilde{\mathcal{M}'}$, from left (dashed purple, shaded gray) to right (solid purple, shaded blue), produces the full meander $\mathcal{M}$.
The required left turns of the inserting arcs, from left to right, lower all Morse numbers again by 1; see \eqref{i}.
This makes the Morse indices of corresponding equilibria coincide, for the solid purple, gray lower left $q$-nest of $\mathcal{M}'$ in (a) and the upper right $q$-nest of $\mathcal{M}$ in (b), respectively.
}}
\label{qnest}
\end{center}
\end{figure}

We start induction on just $r$ with the claim $\mathcal{C}=\mathcal{C}_{rq}$ for $r=1$, in section \ref{r=1}.
That claim has been proved in our prequel \cite{firo23}.
Induction itself is based on two steps.
First we study the suspension $\widetilde{\mathcal{M}'}$ of the meander $\mathcal{M}':=\mathcal{M}^\kappa_{p'q}\,,\ p'=(r-1)(q+1)$; see fig.~\ref{qnest}(a). The second step is an insertion of a $q$-nest, as in fig.~\ref{qnest}(b), to arrive at the target meander $\mathcal{M}:= \mathcal{M}_{pq}\,,\ p=(q+1)$, of our induction.
In sections \ref{susp(r-1)}--\ref{rsubgraph}, we determine the effects of these steps on the associated connection graphs $\mathcal{C}':=\mathcal{C}_{r-1,\,q}$ and $\widetilde{\mathcal{C}'}$ with respect to $\mathcal{C}$; see \eqref{rsub}.
Section \ref{qsubgraph} then utilizes the induction hypothesis on $q-1$.
Equivalence corollary \ref{rqqrcor} swaps the position of $r$ and $q$ to play this case back to the previous sections; see \eqref{connrqqr}--\eqref{qsub}.
Section \ref{tosinks} determines the very few remaining heteroclinic orbits, by the adjacency characterization \eqref{wolfrum}.
This concludes the induction step over $r$ and $q$.

\subsection{The case $r=1$}\label{r=1}
The case $\mathcal{C}^\kappa_{rq}$ with $r=1$ has already been treated in theorem 7.6 of \cite{firo23}.
In terms of the equilibrium labels $A_j'\,,B_j'\,,C_j'\,,D_j'$ used there:
\begin{equation}
\label{connr=1}
\begin{aligned}
A'_j \quad &\leadsto \quad A'_{j-1},B'_{j-1}, \phantom{D'_{j-1},} \quad\textrm{for } 1\leq j \leq q;&&\\
B'_j \quad &\leadsto \quad A'_{j-1},B'_{j-1}, \phantom{D'_{j-1},} \quad\textrm{for } 1\leq j \leq q-1;&&\\
C'_j \quad &\leadsto \quad B'_{j-1},C'_{j-1},D'_{j-1}, \quad\textrm{for } 2\leq j \leq q; \quad &C'_1 \quad &\leadsto \quad B'_0,D'_0;\\
D'_j \quad &\leadsto \quad A'_{j-1},C'_{j-1},D'_{j-1}, \quad\textrm{for } 2\leq j \leq q+1; \quad &D'_1 \quad &\leadsto \quad A'_0,D'_0.\\
\end{aligned}
\end{equation}
See also fig.~7.1 in \cite{firo23}. 
Reverting the trivial equivalence $\kappa$ by $180^\circ$ rotation,
the substitutions 
\begin{equation}
\label{connr=1labels}
\begin{aligned}
A'_0 &\mapsto A^0_1 & B'_0 &\mapsto B^0_0 \qquad & C'_1 &\mapsto B^1_0 \qquad & D'_0 &\mapsto A^1_0 \\[2pt]
A'_j &\mapsto A^0_{j+1} \qquad & B'_j &\mapsto B^0_j & C'_{j+1} &\mapsto B^1_j  & D'_j &\mapsto A^1_j \\[2pt]
A'_q &\mapsto B^0_q &   &                              & D'_{q+1} &\mapsto B^1_q  & D'_q &\mapsto A^1_q \\
\end{aligned}
\end{equation}
for $0< j < q-1$ confirm our claim.
Here we have compared our current $h_1$ labels \eqref{h1labels-oddr}, for $r=1$, with the labels (7.2) of \cite{firo23}, adapted to $h_1^\kappa$ there.

\subsection{Embedding of $\mathcal{C}^\kappa_{r-1,\,q}$ in $\mathcal{C}$, step 1: suspension }\label{susp(r-1)}

From \cite{firo23}, we first recall the precise notion of suspension, for any dissipative meander $\mathcal{M}$. 
See section 3 there, for further illustration and commentary.
We will then apply suspension to the dissipative meander $\mathcal{M}^\kappa_{p',q}$ with $p':=(r-1)(q+1)$, as in fig.~\ref{qnest}(a).

To define suspension, we first label the $N$ vertices $\mathcal{E}=\{E_1,\ldots,E_N\}$ of $\mathcal{M}$, along the horizontal axis, as $E_j:=h_1(j)$. 
For the $N$ ``interior'' vertices $\{\widetilde{E}_1,\ldots,\widetilde{E}_{N}\}$ among all $N+2$ suspension vertices $\widetilde{\mathcal{E}}=\{\widetilde{E}_0,\ldots,\widetilde{E}_{N+1}\}$, we choose a corresponding enumeration $\widetilde{E}_j=\widetilde{h}_1(j+1)$, for $j=1,\ldots,N$.
This embeds old vertices $\mathcal{E}\subset\widetilde{\mathcal{E}}$ into the suspension via the lifting identification
\begin{equation}\label{Elift}
E_j \mapsto \widetilde{E}_j\,.
\end{equation}
Henceforth, we write $E_j=\widetilde{E}_j$ under this identification, for $j=1,\ldots,N$.

We define the \emph{suspension} $\widetilde{\mathcal{M}}$ as an augmentation of $\mathcal{M}$ by two overarching arcs (black in fig.~\ref{qnest}(a)): an upper arc from the first new vertex $\widetilde{E}_0$ to the last old vertex $\widetilde{E}_N=E_N$, and a lower arc from the first old vertex $\widetilde{E}_1=E_1$ to the last new vertex $\widetilde{E}_{N+1}$. 
This extends the previous definition of $\widetilde{h}_1$ to $\widetilde{h}_1(j):=\widetilde{E}_{j-1}$ for $j=1,\ldots,N+2$.
We can now recall proposition 3.1 of \cite{firo23} for the entourage of Sturm permutations $\widetilde{\sigma}$, Morse numbers $i$, zero numbers $z$, and (formal) Morse-adjacent edges of $\widetilde{\mathcal{C}}$ of  $\widetilde{\mathcal{M}}$, as follows.

\begin{prop}\label{suspprop}
For any dissipative, but not necessarily Morse, meander $\mathcal{M}$, the suspension $\widetilde{\mathcal{M}}$ has the following properties, for all $1\leq j,k\leq N,\ j\neq k$: 
\begin{enumerate}[(i)]
  \item $\widetilde{\sigma}(1)=1$ and $\widetilde{\sigma}(N+2)=N+2$;
  \item $\widetilde{\sigma}(j+1) = N+2 - \sigma(j)=\kappa \sigma(j)+1$;
  \item $i(\widetilde{E}_0)=i(\widetilde{E}_{N+1})=0$;
  \item $i(\widetilde{E}_j)=i(E_j)+1$;
  \item $z(\widetilde{E}_j-\widetilde{E}_0) = z(\widetilde{E}_j-\widetilde{E}_{N+1}) = 0$;
  \item $z(\widetilde{E}_j-\widetilde{E}_k) = z(E_j-E_k)+1$;
  \item $\widetilde{E}_j \leadsto \widetilde{E}_k \quad \Longleftrightarrow \quad E_j \leadsto E_k$\,;
  \item $\widetilde{E}_j \leadsto \widetilde{E}_0\,, \widetilde{E}_{N+1}$\,, in case all $i\geq 0$.
\end{enumerate}
\end{prop}

The consequences for Sturm meanders $\mathcal{M}$ are summarized in corollary 3.1 of \cite{firo23}.

\begin{cor}\label{suspcor}
For Sturm meanders $\mathcal{M}$ the following holds true.
\begin{enumerate}[(i)]
\item The suspension $\widetilde{\sigma}\in S_{N+2}$ of any Sturm permutation $\sigma\in S_N$ is Sturm.
\item All $i=1$ equilibria connect heteroclinically, in $\widetilde{\mathcal{C}}$, towards the two polar $i=0$ sinks $\widetilde{E}_0,\widetilde{E}_{N+1}$ in the bottom row.
\item The connection graph $\widetilde{\mathcal{C}}$ of the suspension  $\widetilde{\mathcal{M}}$ contains the connection graph $\mathcal{C}$ of $\mathcal{M}$, lifted to the rows $i\geq 1$.
\end{enumerate}
\end{cor}

With $p':=(r-1)(q+1)$, we apply proposition \ref{suspprop} to the dissipative meander $\mathcal{M}' := \mathcal{M}^\kappa_{p',q}$\,.
Let $\mathcal{C}' := \mathcal{C}^\kappa_{r-1,q}$ denote the connection graph of $\mathcal{M}'$; see fig.~\ref{qnest}(a).
By induction hypothesis \eqref{ind}, we may assume $\mathcal{M}'$ to be Sturm.
Hence proposition \ref{suspprop} and corollary \ref{suspcor} apply to the suspension $\widetilde{\mathcal{M}'}$ of $\mathcal{M}'$, with connection graph $\widetilde{\mathcal{C}'}$.
In particular, the suspension is also Sturm. 
By proposition \ref{suspprop}(iii),(iv), all Morse indices of the suspension are strictly positive, except for the two new sinks $i(A^1_0)=0=i(A^0_1)$.
By corollary  \ref{suspcor}(ii),(iii), the connection graph $\mathcal{C}'$ is embedded in the suspension $\widetilde{\mathcal{C}'}$, at Morse levels raised by 1, and any equilibrium at Morse level 1 possesses edges  (i.e., heteroclinic orbits) to each of the two new sinks:
\begin{equation}
\label{susp1to0}
A^1_1\,, B^1_0\,,A^2_0 \quad \leadsto \quad A^1_0\,, A^0_1 \,.
\end{equation}

\begin{figure}[t!]
\begin{center}
\centering \includegraphics[width=\textwidth]{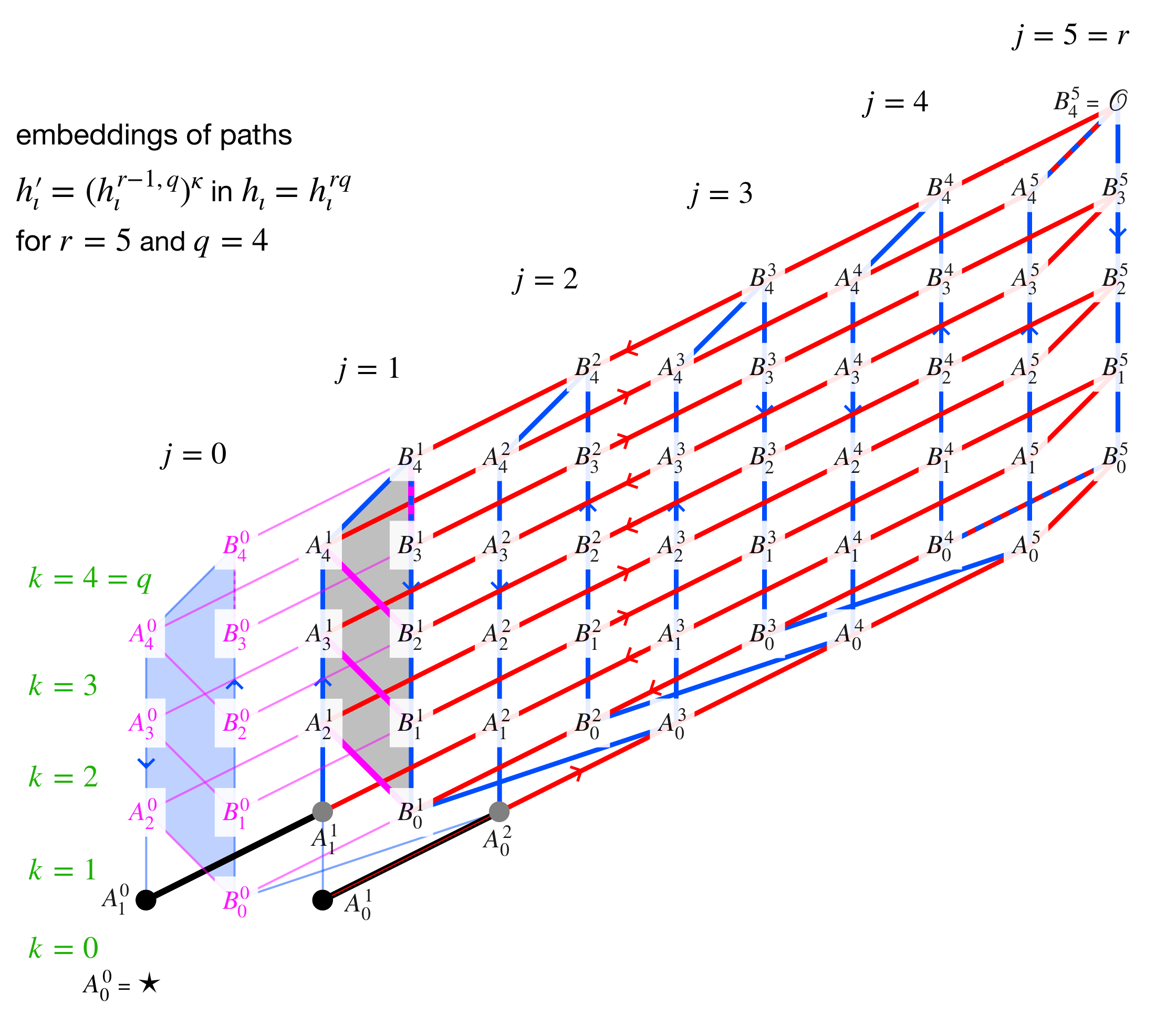}
\caption{\emph{
Induction step, from $(r-1,q)$ to $(r,q)$, on the level of axis and meander paths $h_\iota':=(h^{r-1,q}_\iota)^\kappa$ and $h_\iota:= h^{rq}_\iota\,,\ \iota=0,1$. 
See fig.~\ref{qnest} for the corresponding meanders, and fig.~\ref{h0h1fig} for $h^{rq}_\iota$.
The suspension arcs $A^1_0 A^2_0$ and $A^1_1 A^0_1$ of fig.~\ref{qnest}(a) are colored black again.
The colors, and thicknesses, of $h_0$ and $h_1$ have been modified, compared to fig.~\ref{h0h1fig}.
The $q$-nest insertion in fig.~\ref{qnest}(b) creates the new vertices $B^0_0 \ldots B^0_q A^0_q \ldots A^0_2$ (purple), which bound the vertical Chafee-Infante region $\mathrm{CI}^0$, shaded blue.
\\
The axis path $h_1'$ (thick blue) of meander $\mathcal{M}'$ starts at $A^1_1$ and terminates at $A^2_0$ (both gray).
Thin blue lines mark the extension of $h_1$, beyond the termination of $h_1'$ at $A^2_0$ (gray), towards the termination of $h_1$ at $A^0_1$ (black).
Also note the thin blue single segment of $h_1$ from the start at $A^1_0$ (black) to the subsequent start of $h_1'$ at $A^1_1$ (gray).
On their overlap, the orientations of $h_1'$ and $h_1$ coincide.
\\
The meander path $h_0'$ (thick red) also runs from $A^1_1$ to $A^2_0$.
Due to suspension, however, the orientations of $h_0'$ and $h_0$ are opposite on their overlap; see also fig.~\ref{qnest}(a).
The new upper right $q$-nest of $\mathcal{M}$ in fig.~\ref{qnest}(b) generates a thin purple detour of $h_0$ via $\mathrm{CI}^0$, for each of the original thick purple lower left $q$-nest shortcuts $h_0'$ in $\mathcal{M}'$.
The path $h_0'$ traverses the shortcut arcs $A^1_{j+2} B^1_j$, for $j=0,\ldots,q-2$, and $B^1_q B^1_{q-1}$, in alternating directions.
The shortcut vertices bound the vertical Chafee-Infante region $\mathrm{CI}^1$, shaded gray.
Their $\mathrm{CI}^0$ detours of $h_0$ (thin purple) are via stopover vertices $A^0_{j+2}B^0_j$ and $B^0_qB^0_{q-1}$\,.
}}
\label{embedh0h1}
\end{center}
\end{figure}

We plan to prove claim \eqref{iLjk} on the Morse indices $i$ of $\mathcal{M}$, by induction on $r$.
We therefore need to keep track of $i$, and of the precise vertex labels $L^{rj}_{qk}$\,, during the passages from $\mathcal{M}'$ to $\widetilde{\mathcal{M}}$, and onwards to $\mathcal{M}$.
In fig.~\ref{qnest}(a) we have already abbreviated the equilibrium labels $L^j_k := L^{rj}_{qk}$ according to their appearance on $\mathcal{M}$, in (b).
By induction hypothesis, \eqref{iLjk} holds for the vertices $L^{r-1,j}_{q,k}$ of $\mathcal{M}_{p',q}$\,.
In fig.~\ref{qnest}(a), however, that meander appears as $\mathcal{M}'$, i.e.~rotated by $180^\circ$.
This means that the labels $L^{r-1,j}_{q,k}$ of $\mathcal{M}_{p',q}$ appear in $\mathcal{M}' := \mathcal{M}^\kappa_{p',q}$ in the reverse order of \eqref{h1labels-oddr}, \eqref{h1labels-evenr}, along the horizontal $h_1$ axis: odd superscripts $j$ on the right, and even on the left.
This proves the embedding of vertices by the simple identification
\begin{equation}
\label{Lind}
L^{rj}_{qk}=L^{r-1,\,j-1}_{qk},
\end{equation}
for all $1\leq j\leq r,\ 0\leq k\leq q$ and tags $L=A,B$, except for the suspension vertex $A^1_0$\,.
Since suspension proposition \ref{suspprop}(iv), on the other hand, raises Morse indices by 1, the induction hypothesis implies
\begin{equation}
\label{iLjkind}
i(L^{rj}_{qk})=i(L^{r-1,j-1}_{qk})+1 \ =\ 
\begin{cases}
      \ (j-1)+k-1+1 &=\ j+k-1, \text{ for } L=A, \\
      \ (j-1)+k+1 &=\ j+k, \qquad\, \text{for } L=B,
\end{cases}
\end{equation}
as claimed in \eqref{iLjk}.
Of course, the $i=0$ sinks $A^1_0$ and $A^0_1$ generated by suspension also fit.
In particular, the suspension meander $\widetilde{\mathcal{M}'}$ is Morse, and hence Sturm.
Already, this identifies the correct Morse indices for all vertices in $\mathcal{E}^r_q \setminus \mathrm{CI}^0$.

Much more importantly, however, proposition \ref{suspprop}(vii) also embeds the connection graph $\mathcal{C}_{r-1,q}$ into the suspension $\widetilde{\mathcal{C}'}$, in terms of the labels $L^j_k=L^{rj}_{qk}$ of fig.~\ref{qnest}(a).
Since the labels $L^j_k$ are already borrowed from $\mathcal{C}_{rq}$\,, i.e.the target connection graph of our induction, we describe this embedding in terms of that supposed target $\mathcal{C}$ itself; see \eqref{connclaim} and fig.~\ref{embedconn}.
Indeed, the identification \eqref{Lind} and the induction hypothesis \eqref{ind} imply that the connection graph of the suspension $\widetilde{\mathcal{C}'}$ coincides with $\mathcal{C}$, with the exception of the vertices and edges involving the first vertical Chafee-Infante stack $\mathrm{CIS}^0$. In symbols:\begin{equation}
\label{conntil}
\widetilde{\mathcal{C}'}\setminus \{A^0_1\} \quad = \quad \mathcal{C} \setminus \mathrm{CIS}^0\,.
\end{equation}
Carefully note here that the symbol $\mathcal{C}$ on the right hand side stands for what we \emph{claim} the connection graph $\mathcal{C}_{rq}$ to be; see \eqref{connclaim}.

\subsection{Embedding of $\mathcal{C}^\kappa_{r-1,\,q}$ in $\mathcal{C}$, step 2: insertion of $q$-nest}\label{qnest(r-1)}

We can now insert the lower left $q$-nest of $\widetilde{\mathcal{C}'}$ to form the upper right $q$-nest of the meander $\mathcal{M} := \mathcal{M}_{pq}$\,.
See fig.~\ref{qnest}(b).
Since $p=p'+(q+1)=r(q+1)$, the resulting connection graph $\mathcal{C}$ of $\mathcal{M}$ is the target claim $\mathcal{C}_{rq}=\mathcal{C}$ of our induction.

We aim to extend the Morse indices \eqref{iLjk}, \eqref{iLjkind} to include $\mathrm{CI}^0$.
Trivially, $i(A^0_1)=0$, for all dissipative meanders; see \eqref{i}.
We therefore have to show that, upon $q$-nest insertion,
\begin{enumerate}[(i)]
  \item Morse indices of vertices in embedded $\mathcal{C}'$ are invariant, and
  \item Morse indices of the inserted upper right $q$-nest satisfy \eqref{iLjk}
\end{enumerate}

By fig.~\ref{qnest}, the purple new upper right $q$-nest arcs in (b) arise, by detours in the meander $\mathcal{M}$, from the previous purple shortcut arcs in the lower left $q$-nest of (a), in the meander $\mathcal{M}'$.
Specifically, and ignoring alternating meander orientation, the new arcs $A^0_{j+2} B^0_j$\,, for $j=0,\ldots,q-2$, and the arc $B^0_q B^0_{q-1}$\,, arise from the previous arcs $A^1_{j+2} B^1_j$ and $B^1_q B^1_{q-1}$ by the following detours:
\begin{equation}
\label{detours}
A^1_{j+2}\, A^0_{j+2}\, B^0_j\, B^1_j \quad \text{and} \quad B^1_q\, B^0_q\, B^0_{q-1}\, B^1_{q-1}\,,
\end{equation}
respectively.
See also fig.~\ref{embedh0h1}, where the detours \eqref{detours} in $\mathcal{C}$ are marked thin purple, and the previous shortcuts in $\mathcal{C}'$ are thick purple.
Along each detour, the left turning first arc lowers superscript $j$ and Morse index $i$ by 1.
The following left and right turns lower and raise $i$ by 1, successively.
In particular the Morse indices of the previous vertices remain invariant, as claimed in (i).
Moreover, the Morse indices in the inserted $q$-nest also satisfy \eqref{iLjk}, as claimed in (ii).

As a trivial corollary, we obtain that $A^0_1,\ B^0_0,$ and $A^1_0$ are the only $i=0$ sink vertices of $\mathcal{M}$.


\subsection{Blocking of heteroclinic orbits  in $\mathcal{C}_{rq}$ which emanate from $\mathrm{CI}^0$}
\label{cblock}

The vertices $A^0_1$ and $B^0_0$ are the $i=0$ sinks in $\mathrm{CI}^0$, by the previous section.
No heteroclinic edges emanate from sinks.

All other heteroclinic orbits emanating from $\mathrm{CI}^0$, to any other vertices in $\mathcal{E}^r_q$\,, are blocked by $B^0_0$.
Indeed, evaluation of zero numbers by half rotation numbers, as in \cite{furo91,rofi21}, shows
\begin{equation}
\label{zeval}
z(v-B^0_0) = \begin{cases}
      0, &\text{ for } v\in \mathrm{CI}^0, \\
      0 \text{ or } 1, &\text{ otherwise}.
\end{cases}
\end{equation}
Here we use that $B^0_0$ is the first vertex of $\mathrm{CI}^0$, along the meander path $h_0$, for the first alternative; see \eqref{h0labels-oddq}, \eqref{h0labels-evenq}, and fig.~\ref{embedh0h1}.
Also, $B^0_0$ is the first vertex of $\mathrm{CI}^0$ along the horizontal $h_1$-axis; see \eqref{h1labels-oddr}, \eqref{h1labels-evenr}.
All $v \notin \mathrm{CI}^0$, i.e.~the vertices of the second alternative, precede $B^0_0$ along the horizontal $h_1$-axis.
In particular, $B^0_0$ blocks all heteroclinic orbits from $\mathrm{CI}^0$ to any remaining equilibria $v \notin \mathrm{CI}^0$, by \eqref{zdrop} and \eqref{block}.

\begin{figure}[t!]
\begin{center}
\centering \includegraphics[width=\textwidth]{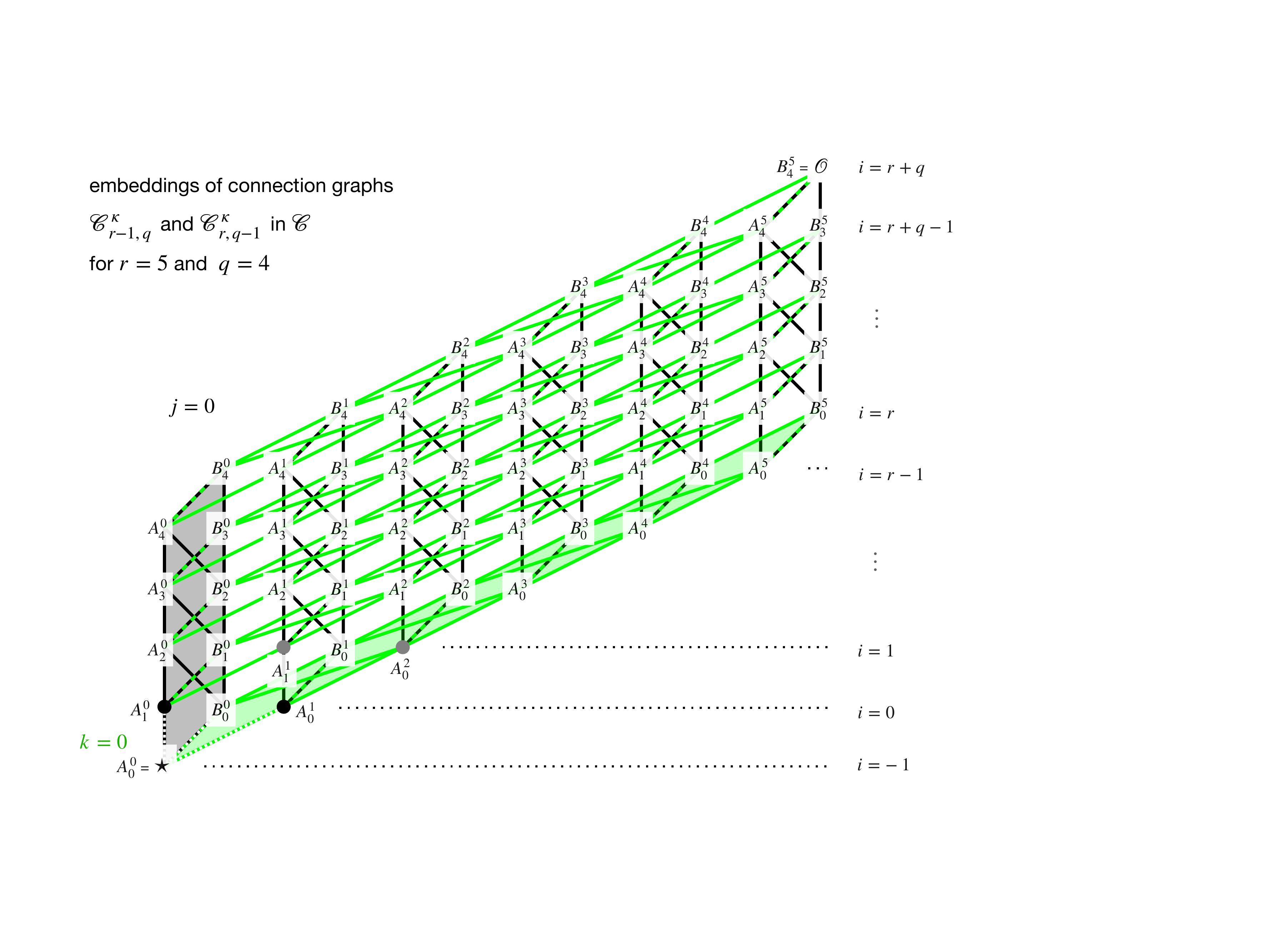}
\caption{\emph{
Induction steps, from $(r-1,q)$ and $(r,q-1)$, to $(r,q)$, on the level of connection graphs $\mathcal{C} := \mathcal{C}_{rq}$\,.
By induction hypothesis on $\mathcal{C}^\kappa_{r-1,q}$\,, and after suspension as in fig.~\ref{qnest}(a), we obtain the connection graph outside the gray shaded region, for Morse levels $i\geq 1$.
In section \ref{rsubgraph} we show that, indeed, the subsequent $q$-nest insertion of fig.~\ref{qnest}(b) does not block any edge in that part $\mathcal{C}^1$ of $\mathcal{C}$.
\\
By induction hypothesis on $\mathcal{C}^\kappa_{r,q-1}$, on the other hand, we obtain the connection graph outside the green shaded region, for Morse levels $i\geq 1$, after suspension as in fig.~\ref{qnest}(a).
Here we invoke the trivial equivalence $\kappa\varrho$ to swap the induction to run from $(q-1,r)$ to $(q,r)$; see lemma \ref{rqqrlem}, corollary \ref{rqqrcor}, and section \ref{qsubgraph}.
The relabeling $\Lambda$ of \eqref{L} then maps the gray shaded region of $(q,r)$ back to the green shaded region of $(r,q)$.
\\
Together, this establishes the connection graph $\mathcal{C}$, for Morse levels $i\geq 1$.
The remaining edges, from Morse level $i=1$ to the sinks $A^0_1, B^0_0$, and $A^1_0$ are established in section \ref{tosinks}.
}}
\label{embedconn}
\end{center}
\end{figure}

\subsection{Embedding of $\mathcal{C}^\kappa_{r-1,\,q}$ in $\mathcal{C}$, step 3: the subgraph}
\label{rsubgraph}

In this section, we show that all heteroclinic orbits of $\widetilde{\mathcal{C}'}$ persist, upon nose insertion, except possibly those towards the termination sink $A^0_1$.
Indeed, Morse numbers in $\widetilde{\mathcal{C}'}$ remain invariant, under $q$-nose insertion, by section \ref{qnest(r-1)}.
Moreover, zero numbers among $\widetilde{\mathcal{C}'}$ remain invariant, because the relevant half-winding numbers of the meander in \cite{furo91,rofi21} remain unaffected by nose insertion.
Finally, the nose insertion occurs to the right of all vertices of $\widetilde{\mathcal{C}'}$, except $A^0_1$; see fig.~\ref{qnest}.
Therefore the inserted $q$-nest cannot block any heteroclinic edge in $\widetilde{\mathcal{C}'} \setminus \{A^0_1\}$\,.
This proves the embedding
\begin{equation}
\label{rsub}
\widetilde{\mathcal{C}^{\,\kappa}_{r-1,\,q}} \setminus \{A^0_1\} \quad = \quad \mathcal{C} \setminus \mathrm{CI}^0.
\end{equation}
As usual, the deletion of the vertices $\mathrm{CI}^0$, on the right, also is meant to delete all edges involving any of those vertices.

\subsection{Embedding of $\mathcal{C}^\kappa_{r,\,q-1}$ in $\mathcal{C}$}
\label{qsubgraph}
So far, we have only embedded the Morse-lifted version of $\mathcal{C}^\kappa_{r-1,\,q}$ in the candidate $\mathcal{C}$ for $\mathcal{C}_{rq}$\,. 
This establishes fig.~\ref{embedconn} outside the vertical gray shaded region of $\mathrm{CI}^0$.
We now address the embedding of the second part of our induction hypothesis \eqref{ind}, i.e.~of $\mathcal{C}^\kappa_{r,\,q-1}$\,.


To start the induction, with $q=1$, we just invoke equivalence corollary \ref{rqqrcor} to observe the trivial equivalence of $\mathcal{C}_{r1}$ with the case $\mathcal{C}_{1r}$ addressed in section \ref{r=1}, albeit for $r$ replaced by the letter $q$ there.

By induction hypothesis \eqref{ind}, the connection graph $\mathcal{C}^\kappa_{r-1,\,q}$ then satisfies theorem \ref{connthm}. 
We will show that this establishes fig.~\ref{embedconn}, and the case of general $r,q$ it stands for, outside the slanted green shaded region of $\mathrm{CI}_0$.
Together with the previous section, this will establish the proper connection graph $\mathcal{C}_{rq}=\mathcal{C}$ at all Morse levels $i \geq 1$.

We invoke the trivial equivalence \eqref{rqqr} of $\sigma_{qr}= \sigma_{rq}^{\varrho\kappa}$, proved in corollary \ref{rqqrcor}, to swap the subscripts $r \leftrightarrow q$ and $r \leftrightarrow q-1$, respectively.
This was achieved via the intertwining relation \eqref{hrqqr} of lemma \ref{rqqrlem}: $\Lambda \,h^{qr}_\iota = h^{rq}_{1-\iota}\, \kappa$.
Here the label map $\Lambda(L^{qk}_{rj}) = L^{rj}_{qk}$ just swaps the pairs of sub- and superscripts, for any vertex tag $L=A,B$.
From corollary \ref{rqqrcor} we recall how $\Lambda$ also induces isomorphisms of connection graphs with arbitrary vertex labels $E_\nu \in \mathcal{E} = \mathcal{E}^q_r$\,:
\begin{equation}
\label{connrqqr}
\begin{aligned}
   (E_1 \leadsto_{qr} E_2) \quad &\Leftrightarrow \quad (\Lambda E_1 \leadsto_{rq} \Lambda E_2)\,, \text{ i.e. }\\[2pt]
\mathcal{C}_{qr} \qquad &\,\cong \qquad \Lambda^{-1}\,\mathcal{C}_{rq}\,,\ \text{ and }\\[2pt]
\mathcal{C}_{q-1,\,r} \qquad &\,\cong \qquad \Lambda^{-1}\,\mathcal{C}_{r,\,q-1}\,.   \\ 
\end{aligned}
\end{equation}
Here the subscripts on the left refer to Sturm permutation $\sigma_{qr}\,,\,\sigma_{q-1,\,r}$\,, whereas on the right they refer to $\sigma_{rq}\,,\,\sigma_{r,\,q-1}$\,.

On the left side of \eqref{connrqqr}, we can therefore recycle the previous arguments of sections \ref{susp(r-1)}--\ref{rsubgraph}, with swapped $r,q$, to obtain the swapped version of embedding \eqref{rsub}:
\begin{equation}
\label{qsubpre}
\widetilde{\mathcal{C}^{\,\kappa}_{q-1,\,r}} \setminus \{A^0_1\} \quad \cong \quad \Lambda^{-1}(\mathcal{C} \setminus \mathrm{CI}_0).
\end{equation}
Applying the label map $\Lambda$ once again, proves the claim
\begin{equation}
\label{qsub}
\widetilde{\mathcal{C}^{\,\kappa}_{r,\,q-1}} \setminus \{A^1_0\} \quad \cong \quad \mathcal{C} \setminus \mathrm{CI}_0\,.
\end{equation}

The connection graph $\mathcal{C}_{rq}$ is associated to the meander $\mathcal{M}_{pq}$\,, which we have have now constructed by (de-)suspension and insertion in two alternative ways.
At all Morse levels $i\geq 1$, the connection graph $\mathcal{C}_{rq}$ is therefore simply given by the union $\mathcal{C}$ of the subgraphs \eqref{rsub} and \eqref{qsub}.
Indeed, the excluded Chafee-Infante sequences $\mathrm{CI}^0$ and $\mathrm{CI}_0$, respectively, only intersect in sink $B^0_0$ and $\star$.
For Morse levels $i\geq 1$, there is no intersection.
Hence the union $\mathcal{C}$ of the subgraphs determines all connectivity within these two sequences, and identifies them as the Chafee-Infante stacks $\mathrm{CIS}^0,\mathrm{CIS}_0$, for $i\geq 1$.
By section \ref{cblock}, no further edges emanate from these stacks.
This proves theorem \ref{connthm}, i.e. $\mathcal{C}_{rq}$=$\mathcal{C}$, at all Morse levels $i\geq$1.

\subsection{Heteroclinic orbits in $\mathcal{C}_{rq}$ from Morse level 1 to 0}
\label{tosinks}
It remains to determine all heteroclinic edges from Morse level $i=1$ to the three $i=0$ sinks $A^0_1,\,B^0_0$, and $A^1_0$ appearing in $\mathcal{C}$;
see fig.~\ref{connfig}.
Adjacency along $h^{rq}_1$ (blue) or $h^{rq}_0$ (red), as in fig.~\ref{h0h1fig}, establishes all claimed edges, except the dashed green/black edges $B^0_1\leadsto A^0_1$ and $B^1_0\leadsto A^1_0$.
These latter edges, however, follow from \eqref{qsub} and \eqref{rsub}, respectively.

To show the absence of any other edges, between $i=1$ and $i=0$, we first recall how section \ref{cblock} blocks all heteroclinic orbits which leave $\mathrm{CI}^0$ or, by section \ref{qsubgraph}, leave $\mathrm{CI}_0$\,.
This proves absence of the four edges $A^2_0A^0_1,\ B^1_0A^0_1,\ B^0_1A^1_0$, and $A^0_2A^1_0$\,.

We prove absence of the fifth and final option $u(t,\cdot): A^1_1 \leadsto B^0_0$ for a heteroclinic edge, by contradiction.
Meander winding and \eqref{z} establish 
\begin{equation}
\label{zA11B00}
 z(A^1_1-B^0_0)=1;
\end{equation}
see also \eqref{zeval}.
Any heteroclinic edge $u(t,\cdot)$, however, would have to satisfy
\begin{equation}
\label{hetA11B00}
1 \ =\ i(A^1_1) \ >\ \lim_{t\rightarrow -\infty} z(A^1_1-u(t,\cdot)) \ \geq 
\ \lim_{t\rightarrow +\infty} z(A^1_1-u(t,\cdot)) \ =\ z(A^1_1-B^0_0) \ =\ 1,
\end{equation}
by dropping \eqref{zdrop} of zero numbers and \eqref{zA11B00}.
This contradiction eliminates any edge from $A^1_1$ to $B^0_0$\,, and completes the induction proof of theorem \ref{connthm}. \hfill $\bowtie$

\section{Discussion: non-Morse cases and continued fractions}\label{Dis}

Let us summarize.
Arc configurations $\mathcal{M}$ are Sturm, if and only if they are dissipative Morse meanders.
The Morse property, in particular, requires nonnegative formal Morse indices, $i_j\geq 0$, in \eqref{i}.
We have established the connection graphs of all primitive $3$-nose meanders $\mathcal{M}_{pq}$ which are dissipative Morse meanders.
The three noses give rise to an upper left $p$-nest, i.e. $p$ nested upper arcs, and a right upper $q$-nest.
All lower arcs are nested, forming a rainbow of $p+q$ lower arcs.
We identify this configuration with the variant of a lower left $q$-nest, a lower right $p$-nest, and an upper rainbow, by trivial equivalence $\kappa$, i.e.~by a $180^\circ$ rotation of the meander.
``Primitive'' means that we do not permit suspension arcs, overarching the upper $p$- and $q$-nests; but see section \ref{susp(r-1)}.
A dissipative meander arises, if and only if $p-1$ and $q+1$ are coprime.
The meander is Morse, if and only if $p=r(q+1)$; see theorem 2.1 in \cite{firo23} and theorem \ref{krmuthm} above.
The (pointed)  connection graph $\mathcal{C}_{rq\star}$ is then given by theorem \ref{connthm}; see also fig.~\ref{connfig}.

Our discussion below aims at dissipative $3$-nose meanders which are not necessarily Morse, nor primitive.
Let $n_0:=p+q$ count the arcs of upper nests.
Since $p-1$ and $q+1$ are coprime, so are $q+1$ and $n_0$\,.
Our approach is based on the continued fraction representation of $n_0/(q+1)$. 
See \cite{per} for background material on continued fractions.
Let $s$ denote the minimal  number of iterated suspensions, which are required to restore the Morse property of the non-Morse permutations, i.e.~to provide Sturm attractors again.
We first review our present results, i.e.~the case of $s=0$ suspensions, from this viewpoint.
We then indicate how trivial isotropy $\sigma^{-1}=\kappa\sigma\kappa$ in the sense of theorem \ref{rqqrthm}, i.e.~involutive $\sigma\kappa$, are related to symmetry of the continued fraction under order reversal.
Time reversibility, in contrast, requires short continued fractions of length three; see \eqref{m=2}.
We calculate the Morse polynomials for short continued fractions, explicitly.
This allows us to identify hopeful non-Morse candidates for time reversibility.
To our surprise, all time-reversible connection graphs which we have encountered by this approach, so far, fall into the same class of Chafee-Infante lattices $\mathcal{C}_{rq\star}$ which we have already encountered in theorem \ref{connthm}.

We do not provide full details on these results, in the present discussion.
Instead, we refer to \cite{firo24}.
To fully expose our continuing ignorance, we conclude with table \ref{n=63} of all $3$-nose Sturm meander configurations with $n=n_0+s=63$ upper arcs and, consequently, with $N=127$ equilibria.
Only three out of the 22 cases listed there are fully covered by our theorems above.

Specifically, we represent the triple $(p,q,n_0)$ by the finite continued fraction
\begin{equation}
\label{cfrac}
\frac{n_0}{q+1} \ =\ b_0 + \frac{1}{b_1 + \frac{1}{b_2 + \frac{1}{\dots + \frac{1}{b_m}}}} \ =\ [b_0,b_1,b_2,\ldots,b_m] \ =:\ b\,.
\end{equation}
Here the integers $b_1,\ldots,b_m$ are strictly positive.
Admitting $b_m=1$ whenever necessary, by convention, we may always assume $m$ to be even.

We claim that the continued fractions \eqref{cfrac}, with $m\geq 2$ and all $b_k\geq 1$, are in 1-1 correspondence with the primitive, but not necessarily Morse, 3-nose meanders.
Indeed, integer $b=[b_0]$ occur for  $n_0=b_0(q+1),\ p = n_0-q = (b_0-1)(q+1)+1$.
Integers $b\geq 2$ do not lead to $3$-nose meanders, because $q+1\geq 2$ then implies that $(p-1,q+1)$ are not coprime.
All Chafee-Infante cases, of only two noses, arise as $b=[0,1,b_2]$, i.e.~via $n_0=q=b_2\geq 1,\ p=0$.
(The case $b=1$ of $n_0=1,\ q=0,\ p=1$ can also be subsumed here, via $b=[0,1,1]$ instead.)
These are the only meanders for which $b_0=0$, i.e. $n_0<q+1$.
This proves our claim

The very special Sturm case $p=r(q+1),\ n_0=(r+1)(q+1)-1$ treated in the present paper just corresponds to
\begin{equation}
\label{rqfrac}
\frac{n_0}{q+1} \ =\ r + \frac{1}{1+1/q} \ =\ [r,1,q]\,.
\end{equation}
By theorem \ref{connthm}, the associated (pointed) connection graph is $\mathcal{C}_{rq\star}$\,.
The Chafee-Infante cases are included, as $p=r=0$.

Let $\sigma_{rq}$ denote the Sturm permutation of \eqref{rqfrac}, as before.
Up to trivial rotation equivalence, the inverse permutation $(\sigma_{rq})^{-1}$ is then given by $\sigma_{qr}\,$, i.e.~by the order reversed continued fraction
\begin{equation}
\label{qrfrac}
\frac{n_0}{(q+1)^*} \ =\ q + \frac{1}{1+1/r} \ =\ [q,1,r]
\end{equation}
with $(q+1)^*=r+1$; see theorem \ref{rqqrthm}.
Here $(q+1)^* \in \{1,\ldots,n_0-1\}$ denotes the inverse of $q+1$, in the multiplicative group $\mathbb{Z}^*_{n_0}$ of elements coprime to $n_0$\,.
Indeed, $(r+1)(q+1)=p+q+1=n_0+1 \equiv 1 \mod n_0$\,, and \eqref{qrfrac} follows explicitly.
The involutive Chafee-Infante case $r=0$ corresponds to the multiplicative unit $1\in \mathbb{Z}^*_{n_0}$\,.

More generally, the continued fractions of multiplicative inverses $a^*\, a \equiv 1 \mod n_0$ are always related by  reversed continued fractions
\begin{equation}
\label{invfrac}
\frac{n_0}{a} \ =\ [b_0,b_1,\ldots,b_m]\quad \Leftrightarrow \quad \frac{n_0}{a^*} \ =\ [b_m,\ldots,b_1,b_0]\,.
\end{equation}
Skipping a few details, \eqref{invfrac} follows from §§9-11 in Perron's classic \cite{per}.
In particular, \eqref{qrfrac} confirms how $3$-nose Sturm permutations $\sigma$ with $\kappa\varrho$-isotropy $\sigma=\sigma^{\kappa\varrho}=\kappa\sigma^{-1}\kappa$, i.e.~with involutive $\sigma\kappa$, arise, if and only if $r=q$.
See also theorem \ref{rqqrthm} and much more cumbersome section \ref{krmupf} above.
The total number $n=n_0$ of arcs in the meander then satisfies
\begin{equation}
\label{r=q}
n+1=(q+1)^2\,.
\end{equation}
See table \ref{n=63} below, for the ``iso'' case $n=n_0=63,\ q=7,\ b=[7,1,7]$ of $N=2n+1=127$ equilibria.

Perhaps more intriguingly, \eqref{qrfrac} generalizes to the trivial equivalence of the permutation $\sigma$ of any primitive $3$-nose meander $\mathcal{M}_{pq}\,$, not necessarily of Morse type, and the ``inverse''  $\sigma^{\kappa\varrho}$.
In other words, the respective sizes of their upper right nests, offset by 1, are multiplicative inverses
\begin{equation}
\label{rhoinv}
(q+1)^*\; (q+1) \equiv 1 \mod n_0\,,
\end{equation}
with correspondingly reversed continued fractions \eqref{cfrac}.
This fact is not quite obvious, and will be proved elsewhere \cite{firo24}.
The isotropic case $\sigma=\sigma^{\kappa\varrho}=\kappa\sigma^{-1}\kappa$, alias $(q+1)^*=q+1$, then arises, if and only if the continued fraction \eqref{cfrac} is symmetric.
Again by \cite{per}, this occurs in general, if and only if \begin{equation}
\label{iso}
n_0 \textrm{ is a divisor of } (q+1)^2-1=q(q+2).
\end{equation}
Again see table \ref{n=63} below for the additional ``iso'' cases $n=63,\ n_0=33,52,60,57$ and $q+1=32,25,41,37$.

In fact, there is a deeper connection between $3$-nose meanders and the continued fraction \eqref{cfrac}.
Roughly speaking, $b_k$ with even index $k$ indicate $b_k$ successive lower right $q$-nest suspension-insertions as in fig.~\ref{qnest} and sections \ref{susp(r-1)}, \ref{qnest(r-1)}.
(Slight modifications are required at $k=0,m$.)
For odd indices $k$, we have to perform $b_k$ correspondingly reflected lower left $p$-nest insertions, each followed by a de-suspension.
The two procedures, by the way, are related to the alternating horizontal and vertical extensions by $b_k$ squares, in the rectangular billiards of \cite{fica13}.
We refer to \cite{firo24}, all the while.

Expectedly, however, de-suspensions give rise to negative ``Morse'' indices, purely formally, by recursion \eqref{i}.
Indeed, the minimal and maximal Morse indices of the primitive $3$-nose meanders $\mathcal{M}_{pq}$ are given by
\begin{equation}
\label{morseminmax}
i_\mathrm{min} = 1-\sum_{k\ \mathrm{odd}}b_k\ ,\qquad i_\mathrm{max} = \sum_{k\ \mathrm{even}}b_k\ ,
\end{equation}
respectively.
Since all $b_k$ are strictly positive, the Morse case $i_\mathrm{min}=0$ reduces to continued fractions \eqref{cfrac} with $m=2$ and $b_1=1$. 
This is just our Sturm case \eqref{rqfrac}.

Proceeding recklessly, i.e.~purely formally, in presence of negative ``Morse'' indices and negative ``zero'' numbers, we may still construct (meaningless) ``connection'' graphs, by the formal blocking rules \eqref{i}--\eqref{block}.
To lend meaning to that purely formal construct, however, we apply the requisite minimum of $s:=-i_\mathrm{min}$ suspensions.
In view of proposition \ref{suspprop}, this defines a Sturm meander $\mathcal{M}^s_{pq}$ with the usual entourage of Sturm permutations, attractors, and pointed connection graphs $\mathcal{C}^s_{pq\star}$\,.
Let
\begin{equation}
\label{morsepol}
\widetilde{M}_\star(x) = \tfrac{1}{x} + \mu_0 + \mu_1 x + \ldots + \mu_d x^d
\end{equation}
denote the associated (pointed) Morse ``polynomial'' of positive degree $d$ with Morse counts $\mu_i\geq0$ and $\mu_d>0$.
By \eqref{morseminmax} and our construction, the degree $d$, which is also the dimension of the global attractor, is simply
\begin{equation}
\label{dim}
d=-1+\sum_k b_k\,. 
\end{equation}
We have to sum over all $k$, this time.
Note how \eqref{dim} agrees with $d=r+q$, in our Sturm case \eqref{rqfrac} of $b_1=1$ with $s=0$ suspensions.
In general, we obtain coefficients
\begin{equation}
\label{morsecounts}
\mu_0=3,\quad \mu_{d-1}=2+\tfrac{1}{2}m,\quad \mu_d=1.
\end{equation}

We are then able to always reconstruct, from the continued fraction $b$, not only $p,q,n_0$ but also, the requisite minimal number $s=-i_\mathrm{min}$ of suspensions to first reach Morse meanders, the total number $n=n_0+s$ of upper (or lower) arcs, and the dimension $d$ of the resulting global attractor.
Indeed, we recall how the coprime numerator $n_0$ and denominator $q+1$ follow from \eqref{cfrac}, as does $p=n_0-q$.
The required suspensions $s$, and therefore $n=n_0+s$, follow from the summation \eqref{morseminmax} of odd entries.
For the dimension $d$, we invoke the total sum \eqref{dim}.

We derive some consequences of \eqref{morsecounts} for our exploration of time reversibility in global attractors of $3$-nose meanders.
Suppose we obtain a time-reversible connection graph by the above construction.
Then the Morse counts must be symmetric, i.e.~$\mu_i=\mu_{d-1-i}$ for all $i$.
In particular, $\mu_0=\mu_{d-1}$ in \eqref{morsecounts} implies $m=2$.
Therefore, time reversibility can only occur for \emph{short continued fractions}, of the form
\begin{equation}
\label{m=2}
\frac{n_0}{q+1}=[b_0,b_1,b_2]\,.
\end{equation}
In passing, we note how $\mu_{d-1}=3$ then justifies our minimal choice of $s=-i_\mathrm{min}$ meander suspensions to encounter reversibility, on any the Sturm level.
Indeed, any larger choice of $s$ would generate only $\mu_0=2$ sink equilibria $i=0$, both created by the last suspension.
See proposition \ref{suspprop}(iii),(iv).
This would violate the necessary condition $\mu_0=\mu_{d-1}$ for time reversibility.

For short continued fractions \eqref{m=2}, i.e.~for $m=2$, the explicit (pointed) Morse polynomial \eqref{morsepol} takes the surprising form
\begin{equation}
\label{morsepol2}
(x-1)^3 \tfrac{x}{x+1}\widetilde{M}_\star(x) \ =\ (x-1)^2(x^{b_0+b_1+b_2}-1)+ x(x^{b_0}-1)(x^{b_1}-1)(x^{b_2}-1)\,.
\end{equation}
As a trivial consequence, for example, we can evaluate the number $n$ of arcs after suspension, directly, as 
\begin{equation}
\label{nb}
n\ =\  \tfrac{1}{2}\widetilde{M}_\star(1)-1 \ =\ b_0+b_1+b_2 + b_0b_1b_2 -1 = d +  b_0b_1b_2\,.
\end{equation}
Indeed, we just divide \eqref{morsepol2} by $(x-1)^3$ and interpret the resulting quotients on the right hand side as derivatives, in the limit $x\rightarrow 1$.
Substituting \eqref{dim} proves the claim.

For the Chafee-Infante case $b=[0,1,b_2]$ with $p=0,\ q=b_2$, we correctly obtain $n=d=b_2$ arcs, and $\widetilde{M}_\star(x)=\tfrac{1}{x}(1+x)(1+\ldots +x^{b_2})$.
For general $b=[b_0,b_1,b_2]$, we obtain
\begin{equation}
\label{m=2pq}
\begin{array}{cll}
\ \ s=b_1-1, &d=b_0+b_1+b_2-1,   &n=d+b_0b_1b_2\,,\\
n_0=n-s, &q=b_1b_2, &p=n_0-q\,.
\end{array}
\end{equation}
Compare also the first nine rows of table \ref{n=63} below.

Note how the Morse polynomial \eqref{morsepol2} is invariant under \emph{any} permutation of the three continued fraction elements $b_0,b_1,b_2$.
(Invariance under just reversal of $[b_0,b_1,b_2]$ would follow from trivial equivalence $\kappa\varrho$, we recall.)
Such invariance encourages us to ask about the connection graphs in case $\{b_0,b_1,b_2\}=\{1,r_1,r_2\}$.
The Morse polynomials then coincide with our time-reversible Sturm class $n_0/(q+1)=[r_1,1,r_2]$ of $p=r_1(r_2+1),\ q=r_2\,,\ n=n_0=r_1r_2+r_1+r_2$; see \eqref{rqfrac}.
Let us therefore address the only related non-Sturm class. 
Up to trivial $\kappa\varrho$-equivalence this is
\begin{equation}
\label{1r1r2}
\frac{n_0}{q+1}=[1,r_1,r_2]
\end{equation} 
with $s=-i_\mathrm{min}=r_1-1$ suspensions; see \eqref{morseminmax}.
Permutation invariance of \eqref{morsepol2} implies $d=r_1+r_2,\ n=r_1r_2+r_1+r_2$, as before.
From \eqref{m=2pq} we obtain $n_0=n-s=r_1r_2+r_2+1,\ q=r_1r_2$, and $p=r_2+1$.

The full Morse polynomial \eqref{morsepol2} factorizes as
\begin{equation}
\label{morse1r1r2}
\begin{aligned}
(x-1)^2\tfrac{x}{x+1}\widetilde{M}_\star(x) \ &=\ (x-1)(x^{1+r_1+r_2}-1)+ x(x^{r_1}-1)(x^{r_2}-1)\ =\\
    &=\ (x^{r_1+1}-1) (x^{r_2+1}-1)\,.
\end{aligned}
\end{equation}
Upon division by $(x-1)^2$, we readily recognize the factorization into the two Chafee-Infante stacks
$(\tfrac{1}{x}+1)(1+\ldots+x^{r_j})$, for $j=1,2$; see also theorem \ref{krmuthm}.
In view of theorem \ref{connthm}, this is not surprising, since the Morse polynomials of the Sturm case $b=[r_1,1,r_2]$ and the properly suspended non-Sturm case $b=[1,r_1,r_2]$ coincide, by permutation invariance of \eqref{morsepol2}.
This tempts us to speculate -- but does not prove in any way -- that the pointed connection graphs of such properly suspended non-Morse meanders might coincide with the time-reversible Chafee-Infante lattices $\mathcal{C}_{r_1r_2\star}$\,.
A few explicit cases look promising, at least.
We hope to settle this question in the foreseeable future \cite{firo24}.

\begin{table}[t!]
\centering
\footnotesize{
\begin{tabular}{||>{$}l<{$}|>{$}c<{$}|>{$}c<{$}||>{$}c<{$}|>{$}c<{$}|>{$}c<{$}|>{$}c<{$}|>{$}c<{$}||>{$}c<{$}|>{$}c<{$}|>{$}r<{$}||}
\hline \hline
 b & p-1 & q+1& d & \textrm{rev} & \textrm{iso} & s & n_0 & (q+1)^* & (p-1)^* & b^* \\ \hline \hline 
[1,1,31]      & 31 & 32	& 32 & \checkmark & -  & 0 & 63 & 2  & 61 & [31,1,1]  \\ 
{[1,31,1]}      & 1	 & 32 & 32  & (\checkmark) & \checkmark & 30 & 33 & 32 & 1 & [1,31,1]  \\ \hline 
[1,3,15]      & 15   & 46 & 18 &  (\checkmark) & - &  2 & 61 & 4 & 57 & [15,3,1]  \\ 
{[1,15,3]}      & 3	 & 46 & 18 &  (\checkmark) & - &  14 & 49 & 16 & 33 & [3,15,1]  \\ 
{[3,1,15]}      & 47 & 16 & 18 &  \checkmark & - &  0 & 63 & 4 & 59 & [15,1,3]  \\ \hline 
[1,7,7]        & 7	 & 50 & 14 &  (\checkmark) & - &  6 & 57 & 8 & 49 & [7,7,1]  \\  
{[7,1,7]}        & 55 & 8 & 14 &  \checkmark & \checkmark &  0 & 63 & 8 & 55 & [7,1,7]  \\ \hline 
[2,2,12]      & 37 & 25 & 15 &  \bcancel{\checkmark} & - &  1 & 62 & 5 & 57 & [12,2,2]  \\ 
{[2,12,2]}      & 27 & 25 & 15 & \bcancel{\checkmark} &  \checkmark & 11 & 52 & 25 & 27 & [2,12,2]  \\ \hline \hline
[1,1,1,8,2]  & 19 & 36 & 12 &  - & - &  8 & 55 & 26 & 29 & [2,8,1,1,1]  \\ \hline
[1,2,1,11,1] & 13 & 38 & 15 &  - & - &  12 & 51 & 47 & 4 & [1,11,1,2,1]  \\ \hline
[3,1,1,1,5]  & 45 & 17 & 10 &  - & - &  1 & 62 & 11 & 51 & [5,1,1,1,3]  \\ \hline
[1,1,1,4,4]  & 21 & 38 & 10 &  - & - &  4 & 59 & 14 & 45 & [4,4,1,1,1]  \\ \hline
[1,4,1,7,1]  & 9	 & 44 & 13 &  - & - &  10 & 53 & 47 & 6 & [1,7,1,4,1]  \\ \hline
[1,2,6,2,1]  & 19 &41  & 11 &  - & \checkmark &  3 & 60 & 41 & 19 & [1,2,6,2,1]  \\ \hline
[1,2,4,1,3]  & 19 & 42 & 10 &  - & - &  2 & 61 & 16 & 45 & [3,1,4,2,1]  \\
{[1,4,2,1,3]}  & 11	 & 48 & 10 &  - & - &  4 & 59 & 16 & 43 & [3,1,2,4,1]  \\ \hline
[3,1,1,2,3]  & 44 & 17 & 9 &  - & - &  2 & 61 & 18 & 43 & [3,2,1,1,3]  \\ \hline
[2,1,3,2,2]  & 39 &22  & 9 &  - & - &  2 & 61 & 25 & 36 & [2,2,3,1,2]  \\  
{[2,2,1,3,2]}  & 34 & 25 & 9 &  - & - &  4 & 59 & 26 & 33 & [2,3,1,2,2]  \\ \hline \hline
[1,1,1,1,5,1,1]  & 24 & 37 & 10 &  - & - & 2 & 61 & 33 & 28 & [1,1,5,1,1,1,1]  \\
{[1,1,1,5,1,1,1]}  & 20 & 37 & 10 &  - & \checkmark & 6 & 57 & 37 & 20 & [1,1,1,5,1,1,1]  \\ \hline \hline
\end{tabular}
}
\caption{\emph{
All $3$-nose Sturm meanders with $n=63$ upper arcs and $127$ equilibria, up to rotation by $180^\circ$.
The column $s$ denotes the number of suspending arcs. 
Columns $p-1$ and $q+1$ indicate the sizes $p$ and $q$ of the upper left and upper right nests, respectively, with a total of $p+q=(p-1)+(q+1)=n_0=n-s$ arcs.
The leftmost column is the continued fraction expansion of $n_0/(q+1)=b=[b_0,\ldots,b_m]$, normalized to even $m$.
The right columns refer to the trivial equivalence which replaces the left Sturm permutation $\sigma$ by $\sigma^{\kappa\varrho}=\kappa\sigma^{-1}\kappa$.
This reverses the continued fraction $b$ to become $b^*$, and produces appropriately modified nest sizes via the multiplicative inverses $(p-1)(p-1)^*\equiv (q+1)(q+1)^*\equiv 1 \mod n_0$\,.
Using $(p-1)\equiv -(q+1) \mod n_0$, we may equivalently invoke $(p-1)(q+1)^*\equiv (q+1)(p-1)^*\equiv -1 \mod n_0$\,.
A checkmark ``$\checkmark$'' in the column ``iso'' indicates the isotropy of symmetric $b=b^*$\,, i.e.~of involutive $\kappa\sigma$.
The attractor dimension $\dim$, as well as the suspension and nest counts $s=-i_\mathrm{min}$ and $n_0$\,, are shared by $b,b^*$ in each row; see \eqref{dim}, \eqref{morseminmax}.
The column ``rev'' checks for Morse reversibility.
Checkmarks ``$\checkmark$'' indicate applications of theorem \ref{connthm}.
Parentheses ``$(\checkmark)$'' indicate case-by-case verification.
The cases where the pointed connection graph \emph{fails to be time-reversible}, even though the pointed Morse polynomial \emph{is reversible}, are indicated by a crossed out checkmark ``$\bcancel{\checkmark}$''.
Note how all $b=[b_0,\ldots,b_m]$ fail to be time-reversible, for $m>2$. 
Moreover, permutations of such $b_k$ produce different Morse polynomials.
Often, they do not even produce the same number $2n+1$ of equilibria, as the absence of most permutations from the table shows.
We have sorted rows by increasing lengths $m+1 \in \{3,5,7\}$ of $b,b^*$\,.
Within each length, we have sorted $b$ lexicographically, up to permutations of $b_k$; lower instantiations on the left.
Horizontal separators have been omitted between permutation related continued fractions $b$.
}}
\label{n=63}
\end{table}

In table \ref{n=63} we illustrate all 3-nose Sturm meanders for the case of $n=63$ upper (or lower) arcs.
Why 63?
We have chosen $n+1=8^2$ to accommodate the ``square'' Morse case $r=q=7$, in view of \eqref{r=q}; see also \eqref{iso}.
Smaller odd squares are squares of primes, and therefore do not accommodate factorizations $n+1=(r_1+1)(r_2+1)$ with $r_1\neq r_2$.
The three smaller even squares did not accommodate all feasible combinations of (ir-)reversibility and (non-)isotropy for the Sturm, non-Sturm, $\min \{b_0,b_1,b_2\}=1,\ \min \{b_0,b_1,b_2\}=2$ and $m\geq4$ cases, as do occur in our table.
This made $n=63$ the smallest ``universal'' choice.

The rows for $q$ are generated starting from the Farey sequence $F_n$ of fractions $n_0/(q+1)$ with coprime numerators $n_0$ and denominators $1<q+1<n_0\leq n$.
The continued fractions $n_0/(q+1)=b=[b_0,\ldots,b_m]$ provide the required minimal numbers of $s=-i_\mathrm{min}$ suspensions, by \eqref{morseminmax}, and the attractor dimension $d$, by \eqref{dim}.
Then the nests involve $n_0=n-s$ arcs.
Alternatively, of course, $n_0$ and $q+1$ follow directly from the Farey fractions, and provide $p-1,\ s$, as well.
The multiplicative inverses $(p-1)^*\,, (q+1)^*$ follow most easily from the order reversed continued fractions $b^*$\,.
Of course they share the same values $d, s, n_0, n$, which are invariant under the trivial equivalence of order reversal.
A checkmark ``$\checkmark$'' in the column ``iso'' indicates the $\kappa\varrho$-isotropic cases of symmetric $b=b^*$, i.e.~of involutive $\sigma\kappa$.

For a more detailed inspection, let us consider length $m+1=3$ of the continued fractions $b,b^*$\,, first.
Then permutation related $b$ share the same Morse polynomial \eqref{morsepol2}.
It is amusing to check the invariance $d+b_0b_1b_2=n=63$ of \eqref{nb}, for the first $9$ rows of table \ref{n=63}, in this context.
See also \eqref{m=2pq}.
Only for the three cases of $b_1=1$, however, theorem \ref{connthm} has established explicit time-reversibility of the connection graphs, i.e.~of the associated Chafee-Infante lattices $\mathcal{C}_{b_0b_2\star}$\,.
This is indicated by a checkmark ``$\checkmark$'' in the column ``rev''.
Reversibility of the remaining four cases with $\min\{b_0,b_1,b_2\}=1$ will be addressed in \cite{firo24}, again via the appropriate Chafee-Infante lattices.
Verification checks ``$(\checkmark)$'' in table \ref{n=63} have been performed manually, on a case-by-case basis, by isomorphism to the pointed connection graphs $\mathcal{C}_{r_1r_2\star}$ of theorem \ref{connthm} present in each such permutation class.

The general case of length $m+1=3$ with all $b_k\geq 2$, at first glance, still looks like a ``time-reversible'' hopeful in terms of the Morse polynomial itself.
Indeed symmetry $\mu_i=\mu_{d-1-i}$ of the Morse counts holds, if and only if
\begin{equation}
\label{morserev}
\widetilde{M}_\star(x)\ =\ x^{d-1}\widetilde{M}_\star(1/x)\,.
\end{equation}
First replacing $x$ by $1/x$ in \eqref{morsepol2}, and then multiplying by $x^{d-1}$ with $d=b_0+b_1+b_2-1$ according to \eqref{dim} easily verifies that claim.
At first, we had therefore placed checkmarks ``$\checkmark$'' in the time reversibility column ``rev'' of table \ref{n=63}, to indicate such Morse reversibility.

However, there are counterexamples to reversibility of the associated connection graphs, in case $\min\{b_0,b_1,b_2\}\geq 2$.
Already the simplest case $b=[2,2,2]$ of $p=8,\ q=4,\ n_0=12,\ s=1,\ d=5, \ n=13$, not included in the table, features \emph{five} $i=1$ saddles with heteroclinic orbits to the left-most sink $A$ of the three sinks at the adjacent Morse level $i=0$.
Any of the three equilibria at the reversed, submaximal Morse level $i=d-1$, in contrast, connects to just \emph{four} equilibria at the lower Morse level $i=d-2$.
This discrepancy of an in-degree \emph{five} at $i(A)=0$, and out-degrees \emph{four} at $i=d-1$, contradicts reversibility of the connection graph.

Consider the simple cases $\{b_0,b_1,b_2\} = \{2,3,4\}$ next. 
By trivial equivalence $\kappa\varrho$ of inverses \eqref{invfrac}, we only have to consider the three cases $b=[2,3,4],\ [2,4,3],\ [3,2,4]$ of $s=2,3,1$ suspensions, respectively.
Comparing in- and out-degrees, again, each of the three connection graphs fails to be time-reversible -- even though their Morse polynomials look alluringly reversible.
The failure is by the same discrepancy between Morse levels $i=0,1$ and $i=d-1,d-2$ as before.
The exact same discrepancies of in- and out-degrees rule out reversibility of the connection graphs for the two remaining permutation related cases $b=[2,2,12]$ and $b=[2,12,2]$ in table \ref{n=63}.
We have therefore indicated these failures by crossing out their two preliminary checkmarks in the column ``rev'': see the entries ``$\bcancel{\checkmark}$''.

Even worse, the connection graphs and their time reversals are all pairwise non-isomorphic as connection graphs, in any of the above sets of permutation related examples.
In particular, none of their time reversed connection graphs can appear from any $3$-nose meander and their suspensions.
So far, we have not found any isomorphisms among connection graphs and their reversals beyond the permutation classes $\min \{b_0,b_1,b_2\}=1$ discussed above, and beyond the trivial equivalences of order reversal in $b$.
All connection graphs associated to reversible classes were of the same type $\mathcal{C}_{r_1r_2\star}$ already encountered in theorem \ref{connthm}, for the representative $b=[r_1,1,r_2]$. Since \eqref{m=2pq} identifies these cases by $n=(r_1+1)(r_2+1)-1$, they simply correspond to the proper divisors of $n+1$.
Table \ref{n=63} clearly illustrates that, for the case $n+1=64=2^6$ and the associated three permutation classes of the divisor pairs $(r_1+1,r_2+1) = (2,32),\ (4,16)$, and $(8,8)$, up to interchanging $r_1$ and $r_2$.

On the pessimistic side, suppose $n+1$ is prime. Then there are no cases of $3$-nose Sturm attractors where our results apply and, probably, reversibility fails in all instances.

For larger lengths $m+1\geq5$, we have already mentioned irreversibility of the Morse polynomials, i.e.~violation of \eqref{morserev} by \eqref{morsecounts}.
But some permutations of the $b_k$, other than just order reversal to $b^*$, even affect the total number $n$ of meander arcs.
This is demonstrated by the marked absence of most such permutation related pairs $b,b^*$ from table \ref{n=63}.
\emph{A forteriori,} such absence demonstrates how even the Morse polynomials do depend on permutations of the $b_k$\,.

At present, therefore, we cannot offer any systematic approach to elucidate the structure of the connection graphs, in that huge maze of remaining open cases.
Indeed, all our results above do not offer more than a first glimpse, so far.
All cases beyond the reversible Chafee-Infante lattices $\mathcal{C}_{r_1r_2\star}$ remain wide open.
And we did not even attempt a \emph{geometric} description of the $3$-nose attractors $\mathcal{A}$, yet, e.g.~in terms of their (signed) Thom-Smale complexes.

Such are the amazing intricacies of Sturm global attractors with just three meander noses, at present.
And we owe it all to the introduction of Sturm meanders in \cite{furo91}, with the help of Giorgio Fusco, more than 30 years ago.



\newpage

\begin{thebibliography}{9999)999}

{\footnotesize{

\bibitem[An86]{an86}
S.~Angenent.
The {M}orse-{S}male property for a semi-linear parabolic equation.
\emph{J.~Diff. Eqns.} \textbf{62} (1986), 427--442.

\bibitem[An88]{an88}
S.~Angenent.
The zero set of a solution of a parabolic equation.
\emph{J. Reine Angew.~Math.} \textbf{390} (1988), 79--96.

%

\bibitem[BrFi89]{brfi89}
P.~Brunovsk\'y and B.~Fiedler.
 Connecting orbits in scalar reaction diffusion equations {II}: The
  complete solution.
 \emph{J.~Diff.~Eqns.} \textbf{81} (1989), 106--135.

\bibitem[ChIn74]{chin74}
N. Chafee and E. Infante.
A bifurcation problem for a nonlinear parabolic equation.
\emph{J. Applic.~Analysis} \textbf{4} (1974), 17--37.


\bibitem[Con78]{con}
C.C.~Conley.
\emph{Isolated Invariant Sets and the Morse Index.} 
CBMS Reg.~Conf.~Ser.~Math. \textbf{38}. AMS, Providence, R.I., 1978.

%

\bibitem[Fi94]{fi94}
B. Fiedler.
Global attractors of one-dimensional parabolic equations: sixteen examples.
\emph{Tatra Mountains Math. Publ.} \textbf{4} (1994), 67--92. 

\bibitem[Fi02]{fi02}
B. Fiedler (ed.).  \emph{{H}andbook of {D}ynamical
{S}ystems} \textbf{2}. {E}lsevier, {A}msterdam 2002.

\bibitem[FiCa13]{fica13} 
B.~Fiedler and P.~Casta\~{n}eda.
Rainbow meanders and Cartesian billiards.
\emph{S\~ao Paulo J.~Math.~Sc.} \textbf{6} (2013), 1--29.

\bibitem[FiGRo14]{fietal14}
B.~Fiedler, C.~Grotta-Ragazzo, and C.~Rocha.
An explicit Lyapunov function for reflection symmetric parabolic differential equations on the circle.
\emph{Russ. Math. Surveys.} \textbf{69} (2014), 419--433.

\bibitem[FiRo96]{firo96}
B. Fiedler and C. Rocha.
Heteroclinic orbits of semilinear parabolic equations. 
\emph{J. Diff. Eqns.} \textbf{125} (1996), 239--281.

\bibitem[FiRo99]{firo99}
B.~Fiedler and C.~Rocha.
Realization of meander permutations by boundary value problems.  {\em
J. Diff. Eqns.} \textbf{156} (1999), 282--308.

\bibitem[FiRo00]{firo00}
B.~Fiedler and C.~Rocha.
Orbit equivalence of global attractors of semilinear parabolic
differential equations.
\emph{Trans. Amer. Math. Soc.} \textbf{352} (2000), 257--284.

%
%

\bibitem[FiRo14]{firo14}
B.~Fiedler and C.~Rocha.
Nonlinear Sturm global attractors: unstable manifold decompositions as regular CW-complexes.
\emph{Discr. Cont. Dyn. Sys.} \textbf{34} (2014), 5099--5122.

\bibitem[FiRo15]{firo15}
B.~Fiedler and C.~Rocha.
Schoenflies spheres as boundaries of bounded unstable manifolds in gradient Sturm systems.
\emph{J.~Dyn.~Diff.~Eqns.} \textbf{27} (2015), 597--626.

\bibitem[FiRo18a]{firo3d-1}
B.~Fiedler and C.~Rocha.
Sturm 3-balls and global attractors 1: Thom-Smale complexes and meanders. 
\emph{S\~ao Paulo J. Math. Sc.} \textbf{12} (2018), 18--67;
doi: 10.1007/s40863-017-0082-8.
arXiv:1611.02003.

\bibitem[FiRo18b]{firo3d-2}
B.~Fiedler and C.~Rocha.
Sturm 3-balls and global attractors 2: Design of Thom-Smale complexes.
\emph{J.~Dyn.~Diff.~Eqns.} (2018); 
doi: 10.1007/s10884-018-9665-z

\bibitem[FiRo18c]{firo3d-3}
B.~Fiedler and C.~Rocha.
Sturm 3-ball global attractors 3: Examples of Thom-Smale complexes. 
\emph{ Discr. Cont. Dyn. Syst. A} \textbf{38} (2018), 3479--3545;
doi: 10.3934/dcds.2018149

\bibitem [FiRo20]{firo20}
B.~Fiedler and C.~Rocha.
Boundary orders and geometry of the signed Thom-Smale complex for Sturm global attractors.
\emph{J.~Dyn.~Diff. Eqns.} (2020); doi: 10.1007/s10884-020-09836-5

\bibitem [FiRo23]{firo23}
B.~Fiedler and C.~Rocha.
Design of Sturm global attractors 1: Meanders with three noses, and reversibility.
\emph{Chaos} \textbf{33}, 083127 (2023); doi: 10.1063/5.0147634

\bibitem [FiRo24]{firo24}
B.~Fiedler and C.~Rocha.
Design of Sturm global attractors 3: Negative Morse indices, suspensions, and time-reversibility of connection graphs.
\emph{In preparation} (2024).


\bibitem[FuRo91]{furo91}
G.~Fusco and C.~Rocha.
 A permutation related to the dynamics of a scalar parabolic {PDE}.
 \emph{J. Diff. Eqns.} \textbf{91} (1991), 75--94.



%
%

\bibitem[He85]{he85}
D.~Henry.
Some infinite dimensional {M}orse-{S}male systems defined by parabolic differential equations.
 \emph{J. Diff. Eqns.} \textbf{59} (1985), 165--205.

\bibitem[Hu11]{hu11}
B.~Hu.
\emph{Blow-up Theories for Semilinear Parabolic Equations.}
Lect. Notes Math. \textbf{2018}, Springer-Verlag, Berlin 2011.

\bibitem[LaBe22]{labe22}
Ph.~Lappicy and Ester Beatriz.
An energy formula for fully nonlinear degenerate parabolic equations in one spatial dimension.
\emph{arXiv:} 2201.04215 (2022).

\bibitem[LaFi19]{lafi18}
P.~Lappicy and B.~Fiedler.
 A Lyapunov function for fully nonlinear parabolic equations in one spatial variable. 
São Paulo J.~Math. Sc. \textbf{13} (2019), 283--291; doi: 10.1007/s40863-018-00115-2
 
\bibitem[Ma78]{ma78}
H.~Matano.
Convergence of solutions of one-dimensional semilinear parabolic equations.
\emph{J. Math. Kyoto Univ.} \textbf{18} (1978), 221--227.

\bibitem[Ma82]{ma82}
H.~Matano.
Nonincrease of the lap-number of a solution for a one-dimensional
semi-linear parabolic equation.
\emph{J. Fac. Sci. Univ. Tokyo Sec. IA} \textbf{29} (1982), 401--441.

\bibitem[MaNa97]{mana97}
H.~Matano and K.-I.~Nakamura.
The global attractor of semilinear parabolic equations on ${S^1}$.
\emph{Discr. Cont. Dyn. Sys.} \textbf{3} (1997), 1--24.

\bibitem[MisMro02]{mis}
K.~Mischaikow and M.~Mrozek.
Conley Index. 
In \cite{fi02}, 393--460.

\bibitem[Pe77]{per}
O.~Perron.
\emph{Die Lehre von den Kettenbrüchen. Band I: Elementare Kettenbrüche.}
Springer, Wiesbaden 1977.


\bibitem[RoFi21]{rofi21} 
C.~Rocha and B.~Fiedler. 
Meanders, zero numbers and the cell structure of Sturm global attractors.
\emph{J. Dyn. Diff. Eqns.} (2021); doi: 10.1007/s10884-021-10053-x

\bibitem[St1836]{st1836}
C.~Sturm.
Sur une classe d'\'equations \`a diff\'erences partielles.
\emph{J.~Math.~Pure~Appl.} \textbf{1} (1836), 373--444.
 
 \bibitem[Wo02]{wo02}
M.~Wolfrum.
Geometry of heteroclinic cascades in scalar parabolic differential equations.
\emph{J. Dyn. Diff. Eqns.} \textbf{14} (2002), 207--241.

\bibitem[Ze68]{ze68}
T.I. Zelenyak.
Stabilization of solutions of boundary value problems for a second order parabolic equation with one space variable.
\emph{Diff.~Eqns.} \textbf{4} (1968), 17--22.

\bibitem[zb23]{zb}
Zentralblatt MATH, zbmath.org.
Subject classification MSC 35K57, July 2023.

}}


\end{thebibliography}
\end{document}